\documentclass[11pt,a4paper,oneside]{amsart}
\makeatletter
\@namedef{subjclassname@2010}{%
  \textup{2010} Mathematics Subject Classification}
\makeatother
\usepackage{a4wide}
\usepackage{enumitem}
\usepackage{amsfonts}
\usepackage{amssymb}
\usepackage{amsmath}
\usepackage{amsthm}
\usepackage[all]{xy}
\usepackage{graphicx}
\usepackage{mathrsfs}
\usepackage{color}
\usepackage{url}
\usepackage{array}
\usepackage{pifont}
\usepackage{tikz}
\usepackage{caption}

\allowdisplaybreaks

\newtheorem{thm}{Theorem}[section]
\newtheorem{lem}[thm]{Lemma}
\newtheorem{cor}[thm]{Corollary}
\newtheorem{note}[thm]{Note}
\newtheorem{mukaispecseq}[thm]{Mukai Spectral Sequence}
\newtheorem{dualspecseq}[thm]{``Duality'' Spectral Sequence}
\newtheorem{specseq}[thm]{Spectral Sequence}

\newtheorem{prop}[thm]{Proposition}
\newtheorem{conj}[thm]{Conjecture}

\newtheorem{defn}[thm]{Definition}
\newtheorem{notation}[thm]{Notation}
\theoremstyle{remark}
\newtheorem{rem}[thm]{Remark}
\newtheorem{exmp}[thm]{Example}
\newtheorem*{theorem*}{Theorem}
\newtheorem*{notation*}{Notation}
\theoremstyle{remark}

\numberwithin{equation}{section}

\newcommand{\coh}{\operatorname{Coh}}
\newcommand{\Supp}{\operatorname{Supp}}

\newcommand{\ch}{\operatorname{ch}}

\newcommand{\Tor}{\operatorname{Tor}}
\newcommand{\Ext}{\operatorname{Ext}}

\newcommand{\calHom}{\operatorname{\mathcal{H}\textit{om}}}
\newcommand{\calExt}{\operatorname{\mathcal{E}\textit{xt}}}
\newcommand{\Hom}{\operatorname{Hom}}
\newcommand{\coker}{\operatorname{coker}}
\newcommand{\im}{\operatorname{im}}
\newcommand{\id}{\operatorname{id}}
\newcommand{\rk}{\operatorname{rk}}
\newcommand{\RR}{\operatorname{\mathbf{R}}}
\newcommand{\LL}{\operatorname{\mathbf{L}}}
\newcommand{\Pic}{\operatorname{\mathbf{Pic}}}
\newcommand{\HN}{\operatorname{HN}}
\newcommand{\NS}{\operatorname{NS}}

\newcommand{\HPsi}{\widehat{\Psi}}

\newcommand{\om}{\omega}

\newcommand{\A}{\mathcal{A}}
\newcommand{\B}{\mathcal{B}}

\newcommand{\F}{\mathcal{F}}
\newcommand{\I}{\mathcal{I}}
\newcommand{\OO}{\mathcal{O}}

\newcommand{\T}{\mathcal{T}}

\newcommand{\SC}{\mathscr{C}}

\newcommand{\SM}{\mathscr{M}}

\newcommand{\Po}{\mathscr{P}}

\newcommand{\C}{\mathbb{C}}
\newcommand{\N}{\mathbb{N}}
\newcommand{\Q}{\mathbb{Q}}
\newcommand{\R}{\mathbb{R}}

\newcommand{\Z}{\mathbb{Z}}

\begin{document}

\title[FMTs and Stability Conditions on  Abelian Threefolds]
{Fourier-Mukai Transforms and Bridgeland Stability Conditions on Abelian Threefolds}
\author[Antony Maciocia]{Antony Maciocia}
\address{School of Mathematics\\University of Edinburgh\\
The King's Buildings\\Mayfield Road\\Edinburgh\\EH9 3JZ\\UK.}
\email{A.Maciocia@ed.ac.uk}
\author[Dulip Piyaratne]{Dulip Piyaratne}
\address{Kavli Institute for the Physics and Mathematics of the Universe\\
The University of Tokyo\\
5-1-5 Kashiwanoha\\
Kashiwa\\ 277-8583\\ Japan.}
\email{dulip.piyaratne@ipmu.jp} 

\date{\today}



\begin{abstract}
We use the ideas of Bayer, Bertram, Macr\'i and Toda to construct a Bridgeland stability condition
on a principally polarized abelian threefold $(X,L)$ with $\NS(X) =
\Z[\ell]$ by establishing their Bogomolov-Gieseker type
inequality for certain tilt stable objects associated to the
pair $(\A_{\frac{\sqrt{3} }{2}\ell, \frac{1}{2}\ell}, Z_{\frac{\sqrt{3}
  }{2}\ell, \frac{1}{2}\ell})$ on $X$. This is done by proving the
stronger result that  $\A_{\frac{\sqrt{3}}{2}\ell,\frac{1}{2}\ell}$ is
preserved by a suitable Fourier-Mukai transform.
\end{abstract}

\maketitle

\section*{Introduction}
\label{section0}
In \cite{Bri1} Bridgeland  introduced the notion of stability
conditions on triangulated categories and these now have many
applications to the study of the geometry of the underlying spaces
and highlight the role played by the derived categories of the
suitable categories of sheaves on the spaces.  The space of stability conditions is known precisely for
curves and for abelian surfaces and Bridgeland's geometric stability
conditions provide examples for all projective surfaces (see, for example, \cite{Bri2}, \cite{Macri},
\cite{ABL}, \cite{Okada}).
 A conjectural construction of
Bridgeland stability conditions for projective threefolds was introduced in \cite{BMT} and the
problem is reduced to proving an inequality, which the authors call a Bogomolov-Gieseker
(B-G for short) type inequality, holds for certain tilt stable
objects. This inequality has been shown to hold for three dimensional projective
space (see \cite{BMT} and \cite{Macri1}) and 
smooth quadric threefold (see \cite{Sch}), and 
some progress has been made for more general threefolds (see \cite{Tod1} and \cite{LM}).
 However, there is no known example of a stability condition
on a projective Calabi-Yau threefold and this case is especially
significant because of the interest from Mathematical Physics and also
in connection with Donaldson-Thomas invariants.
In this paper, we establish the existence of a particular stability
condition on a particular Calabi-Yau threefold (namely, a principally
polarized abelian threefold with Picard rank one). However, it is
likely that the method will generalize to other Calabi-Yau threefolds
while the extension to other stability conditions for the abelian
threefold case will be the subject of a forthcoming article.

We reduce the requirement of the B-G type inequality to
a smaller class of tilt stable objects as defined in the Definition
\ref{note2.0}. Moreover, they are essentially minimal objects (also called
simple objects in the literature) of the heart of the stability condition. In
this paper we use Fourier-Mukai theory to prove the B-G
type inequality for these minimal objects by showing that the heart is
preserved by a suitable Fourier-Mukai transform (or FMT for short). For the surface case,
the fact that a countable family of (Bridgeland's) geometric
stability conditions satisfy the numerical conditions for being a stability
condition is actually equivalent to the existence of a Fourier-Mukai
transform preserving the heart. The forward implication was proved by
Huybrechts (\cite{Huy2}) and the reverse implication is a fairly
straightforward exercise (partly done in \cite{Yosh1}). For the
threefold case, we build on these ideas to establish the reverse
implication for our case.

Throughout this paper our abelian varieties will be principally
polarized abelian varieties with Picard rank one over $\C$. Let $(X,L)$ be
an abelian variety of dimension three and let $\ell$ be $c_1(L)$. We
use $L$ to canonically identify $X$ with $\Pic^0(X)$. Let $\Phi : D^b(X) \to
D^b(X)$ be the (classical) FMT with the Poincar\'{e} line bundle on $X \times X$
as the kernel. Then the image of the category $\coh(X)$ under the FMT
$\Phi$ is a subcategory of $D^b(X)$ with non-zero
$\coh(X)$-cohomologies in $0,1,2$ and $3$ positions.
In section 4, we study the slope stability of $\coh(X)$-cohomologies under the
FM transform $\Phi$ in great detail. In particular, we investigate the
images under $\Phi$ of torsion sheaves supported in 
dimensions $1$ and $2$, and of torsion free sheaves whose
Harder-Narasimhan (or H-N for short) semistable factors satisfy
certain slope bounds.

In \cite{BMT} and \cite{BBMT}, the authors construct their conjectural stability
condition hearts as a tilt of a tilt. The first tilt of $\coh(X)$
associated to the H-N filtration with respect to the twisted slope
$\mu_{\omega,B}$ stability is denoted
$\B_{\omega,B}$ and the second $\A_{\omega,B}$ associated to the
H-N filtration with respect to the tilt slope
$\nu_{\omega,B}$ stability. We shall consider the
particular case where $\omega=\sqrt{3}\ell/2$ and $B=\ell/2$.
Let $\Psi:
= L \Phi$ and $\HPsi := \Phi L^{-1}[1]$. At the end of section 4 we prove the images of
the abelian category $\B_{\frac{\sqrt{3}}{2}\ell, \frac{1}{2}\ell}$
under the Fourier-Mukai transforms $\Psi$ and $\HPsi$ have non-zero
$\B_{\frac{\sqrt{3}}{2}\ell, \frac{1}{2}\ell}$-cohomologies only in
positions $0$, $1$ and $2$ (see Theorem \ref{prop4.17}).
On the other hand, we have the isomorphisms (\cite{Muk2})
\begin{align*}
\Psi \circ \HPsi \cong (-1)^* \id_{D^b(X)}[-2], \text{ and }  \HPsi \circ \Psi \cong (-1)^* \id_{D^b(X)}[-2].
\end{align*}
Therefore the abelian category $\B_{\frac{\sqrt{3}}{2}\ell,
  \frac{1}{2}\ell}$ behaves somewhat similarly to the category of
coherent sheaves on an abelian surface under the Fourier-Mukai
transform (see \cite{BBR}, \cite{Macio1}, \cite{Yosh1} for further
details).
Now Theorem \ref{prop4.17} becomes the key technical tool to show that the second tilt
$\A_{\frac{\sqrt{3}}{2}\ell, \frac{1}{2}\ell}$ is preserved by $\Psi$.

Under this auto-equivalence, minimal objects are mapped to minimal
objects and this provides us with an inequality which bounds the top
component of the Chern character of the object. This is the main idea to show that the B-G
type inequality is satisfied by our restricted class of minimal objects in
$\A_{\frac{\sqrt{3}}{2}\ell, \frac{1}{2}\ell}$. In section 5, we have
to show that the B-G type inequality is satisfied by a very special
class of minimal objects by showing that they actually do not exist. This
result is of interest in its own right as it shows that if a bundle
$E$ of such a threefold satisfies $c_1(E)=0=c_2(E)$ then it cannot
carry a non-flat Hermitian-Einstein connection.

\section*{Notation}
\begin{enumerate}[label=(\roman*)]
\item For $0 \le i \le \dim X$ \begin{align*}
\coh^{\le i}(X) &: = \{E \in \coh(X): \dim \Supp(E) \le i  \},\\
\coh^{\ge i}(X) &:= \{E \in \coh(X): \text{for } 0 \ne F \subset E,   \ \dim \Supp(F) \ge i  \}, \text{ and} \\
\coh^{i}(X) &: = \coh^{\le i}(X) \cap \coh^{\ge i}(X).
\end{align*}

\item For an interval $I \subset \R\cup\{+\infty\}$,
$\HN^{\mu}_{\om, B}(I) := \{ E \in \coh(X) : [\mu_{\om , B}^{-}(E) , \mu_{\om , B}^{+}(E)] \subset I \}$.
Similarly the subcategory $\HN^{\nu}_{\om, B}(I) \subset \B_{\om,B}$
is defined.

\item For a FMT $\Upsilon$ and a heart $\mathfrak{A}$ of a t-structure
  for which
$D^b(X)\cong D^b(\mathfrak{A})$, $\Upsilon^k_{\mathfrak{A}}(E) :=
  H^{k}_{\mathfrak{A}}(\Upsilon(E))$.

\item For a sequence of integers $ i_1, \ldots, i_s$,
\[V^{\Upsilon}_{\mathfrak{A}}(i_1, \ldots, i_s) := \{ E \in D^b(X) :
  \Upsilon^j_{\mathfrak{A}}(E) = 0 \text{ for } j \notin \{i_1,
  \ldots, i_s \} \}.
\]
Then $E \in \coh(X)$ being $\Upsilon$-WIT$_i$ is equivalent
  to $E \in V^{\Upsilon}_{\coh(X)}(i)$.

\item Let $(X, L)$ be a principally polarized abelian variety. Then we write
  $\Phi$ for the FMT from $X$ to $X$ with the Poincar\'{e} line bundle
  $\Po:= m^*L \otimes p_1^*L^{-1} \otimes p_2^*L^{-1}$ on $X \times X$
  as the kernel.

\item For $E \in \coh(X)$, $E^k := \Phi^k_{\coh(X)}(E)$.

\item $\Psi: = L \Phi$ and $\HPsi := \Phi L^{-1}[1]$. Here and elsewhere we abuse
  notation to write $L$ for the functor $L\otimes {-}$.

\item For  a  polarized projective threefold $(X,L)$ with Picard rank
  1 over $\C$, the Chern character of $E$ is $\ch(E) = (a_0, a_1 \ell,
  a_2 \frac{\ell^2}{2}, a_3 \frac{\ell^3}{6})$ for some $a_i \in
  \Q$. For simplicity we write $\ch(E) = (a_0, a_1, a_2, a_3)$. Here
  $a_i \in \Z$ for the principally polarized abelian threefold case.
\end{enumerate}
\section{Preliminaries}
\label{section1}
\subsection{Construction of stability conditions}
We recall the conjectural construction of stability conditions as introduced in \cite{BMT}.

Let $X$ be a smooth projective  threefold over $\C$. Let $\om, B$ be in
$\NS_{\R}(X)$ with $\om$ an ample class. The twisted Chern character
$\ch^B$ with respect to $B$ is defined by $\ch^B(-) = e^{-B} \ch
(-)$. So we have
\begin{align*}
& \ch^B_0 =  \ch_0,&    &\ch^B_1 =  \ch_1 - B \ch_0, \\
& \ch^B_2 =  \ch_2 - B \ch_1 + \frac{B^2}{2}\ch_0,&  &\ch^B_3 =  \ch_3
- B \ch_2 + \frac{B^2}{2}\ch_1 -\frac{B^3}{6}\ch_0.
\end{align*}
The twisted slope $\mu_{\om , B}$ on $\coh(X)$ is defined by
$$
\mu_{\om , B} (E) =
\begin{cases}
+ \infty & \text{if } E \text{ is a torsion sheaf} \\
\frac{\om^{2} \ch_1^B(E)}{\ch^B_0(E)} & \text{otherwise}
\end{cases}
$$
for $E \in \coh(X)$.
Then $E$ is said to be $\mu_{\om , B}$-(semi)stable, if for any $0 \ne
F \varsubsetneq E$, we have $\mu_{\om , B}(F)< (\le) \mu_{\om ,
  B}(E/F)$. The H-N filtration of $E$ with respect to $\mu_{\om ,
  B}$-stability enables us to define the following slopes:
\begin{align*}
\mu_{\om , B}^{+}(E)  = \max_{0 \ne G \subseteq E} \ \mu_{\om , B}(G), \ \ \ \
\mu_{\om , B}^{-}(E)  = \min_{G \subsetneq E} \ \mu_{\om , B}(E/G).
\end{align*}
For an interval $I \subset \R\cup\{+\infty\}$, the subcategory $\HN^{\mu}_{\om, B}(I) \subset \coh(X)$ is defined by
$$
\HN^{\mu}_{\om, B}(I) = \{ E \in \coh(X) : [\mu_{\om , B}^{-}(E) , \mu_{\om , B}^{+}(E)] \subset I \}.
$$
Define the subcategories  $\T_{\om , B}$ and $\F_{\om , B}$ of $\coh(X)$ by setting
\begin{align*}
\T_{\om , B} = \HN^{\mu}_{\om, B}(0, +\infty], \ \ \
\F_{\om , B} = \HN^{\mu}_{\om, B}(-\infty, 0].
\end{align*}
Then  $( \T_{\om , B} , \F_{\om , B})$ forms a torsion pair on $\coh(X)$.
Let the abelian category $\B_{\om , B} = \langle \F_{\om , B}[1] ,
\T_{\om , B} \rangle \subset D^b(X)$ be the corresponding tilt of
$\coh(X)$.

The central charge function $Z_{\om , B} : K(X) \to \C$ is defined by
$$
Z_{\om , B}(E) = - \int_X e^{-B - i \omega} \ch(E).
$$
So  $Z_{\om , B}(E) = \left(- \ch^B_3 (E) + \frac{\omega^2}{2}
  \ch^B_1(E) \right) + i \left( \omega\ch^B_2(E)- \frac{\omega^3}{6}
  \ch^B_0(E) \right)$.
The following result is very useful:
\begin{lem} \cite[Lemma 3.2.1]{BMT}
\label{prop1.1}
For any $0 \ne E \in \B_{\om, B}$, one of the following conditions holds:
\begin{enumerate}[label=(\roman*)]
\item $\om^2 \ch_1^B(E) > 0$,
\item $\om^2 \ch_1^B(E) =0$ and $\Im \,Z_{\om, B}(E) > 0$,
\item  $\om^2 \ch_1^B(E) = \Im \,Z_{\om, B}(E) =0$, $- \Re \,Z_{\om,
    B}(E) > 0$ and $E \cong T$ for some $0 \ne T \in \coh^0(X)$.
\end{enumerate}
\end{lem}
As a result of this Lemma, they go on to remark that the vector $(\om^2 \ch_1^B, \Im\,Z_{\om,
  B}, - \Re\,Z_{\om, B})$ for objects in $\B_{\om, B}$ behaves like
the Chern character vector $\ch = (\ch_0, \ch_1 , \ch_2)$ for coherent
sheaves on a surface.

Following \cite{BMT}, the tilt-slope $\nu_{\om , B}$ on $\B_{\om , B}$ is defined by
$$
\nu_{\om , B}(E) =
\begin{cases}
+\infty & \text{if } \omega^2 \ch^B_1(E) = 0 \\
\frac{\Im\,Z_{\om , B} (E)}{\omega^2 \ch^B_1(E)} & \text{otherwise}
\end{cases}
$$
for $E \in \B_{\om , B}$.
Then $E$ is said to be $\nu_{\om , B}$-(semi)stable, if for any $0 \ne
F \varsubsetneq E$ in $\B_{\om , B}$, we have $\nu_{\om , B}(F)< (\le)
\nu_{\om , B}(E/F)$.
In \cite{BMT} it is proved that the abelian category $\B_{\om , B}$
satisfies the H-N property with respect to the tilt-slope
stability. So the following slopes can be defined for  $E \in \B_{\om
  , B}$:
\begin{align*}
\nu^{+}_{\om , B}(E)  = \max_{0 \ne G \subseteq E} \ \nu_{\om , B}(G), \ \ \ \
\nu^{-}_{\om , B}(E)  = \min_{G \subsetneq E} \ \nu_{\om , B}(E/G).
\end{align*}
For an interval $I \subset \R\cup\{+\infty\}$, the subcategory
$\HN^{\nu}_{\om, B}(I) \subset \B_{\om, B}$ is defined by
$$
\HN^{\nu}_{\om, B}(I) = \{ E \in \B_{\om, B} : [\nu_{\om , B}^{-}(E) , \nu_{\om , B}^{+}(E)] \subset I \}.
$$
Define the subcategories $\T_{\om , B}'$ and $\F_{\om , B}'$ of $\B_{\om, B}$ by setting
\begin{align*}
\T_{\om , B}' = \HN^{\nu}_{\om, B}(0, +\infty], \ \ \
\F_{\om , B}' = \HN^{\nu}_{\om, B}(-\infty, 0].
\end{align*}
Then $( \T_{\om , B}' , \F_{\om , B}')$ forms a torsion pair on $\B_{\om , B}$.
Let the abelian category $\A_{\om , B} = \langle \F_{\om , B}'[1] ,
\T_{\om , B}' \rangle \subset D^b(X)$ be the corresponding tilt of
$\B_{\om , B}$.

\begin{conj} \cite[Conjecture 3.2.6]{BMT}
\label{prop1.2}
The pair $(Z_{\om , B}, \A_{\om , B} )$ is a Bridgeland stability condition on $D^b(X)$.
\end{conj}
\begin{defn}
Let $\SC_{\om , B}$ be the class of $\nu_{\om , B}$-stable objects  $E \in \B_{\om , B}$ with $\nu_{\om , B}(E) = 0$.
\end{defn}
Then $E[1] \in \A_{\om , B}$ for any $E \in \SC_{\om , B}$.
\begin{conj} \cite[Conjecture 3.2.7]{BMT}
\label{prop1.3}
Any $E \in \SC_{\om , B}$ satisfies the so called \textbf{Bogomolov-Gieseker Type Inequality}:
$$
\Re\,Z_{\om , B} (E[1]) < 0, \text{ i.e. } \ch^B_3(E) <  \frac{\om^2}{2} \ch^B_1(E).
$$
\end{conj}
Assume $B \in \NS_{\Q}(X)$ and  $\om \in \NS_{\R}(X)$ be an ample class with $\om^2$ is rational. Then the abelian category $\A_{\om, B}$ satisfies the following important property. This was proved for rational classes $\om$ in \cite{BMT}. However a similar proof can be used when we have a weaker condition, namely $\om^2$ is rational. For example, a different parametrization given by $\om \mapsto \sqrt{3} \om$ is considered in \cite{Macri1}.
\begin{lem} \cite[Proposition 5.2.2]{BMT}
\label{prop1.4}
The abelian category $\A_{\om, B}$ is Noetherian.
\end{lem}
As a corollary we have the following
\begin{cor} \cite[Corollary 5.2.4]{BMT}
\label{prop1.5}
The Conjectures \ref{prop1.2} and \ref{prop1.3} are equivalent.
\end{cor}
\subsection{Fourier-Mukai transforms on abelian threefolds}
Let us quickly recall the notion of Fourier-Mukai transform on abelian threefolds. See \cite{BBR}, \cite{Huy1} for further details on Fourier-Mukai theory.

Let $(X, L)$ be a principally polarized abelian threefold with Picard rank 1. Let $\ell := c_1(L)$. Then $\chi(L) = \frac{\ell^3}{6} =1$.
Let $\Po = m^*L \otimes p_1^* L^{-1} \otimes p_2^* L^{-1}$ be the Poincar\'{e} line bundle on $X \times X$.
Then the Fourier-Mukai transform $\Phi : D^b(X) \to D^b(X)$ with kernel $\Po$ is defined by
$$
\Phi(-) := \RR p_{2*} (\Po \stackrel{\LL}{\otimes} p_1^*(-)).
$$
Here $X \stackrel{p_1}{\longleftarrow}  X \times X \stackrel{p_2}{\longrightarrow}  X $ are the projection maps.
In \cite{Muk2} Mukai proved that $\Phi$ is an auto-equivalence of the derived category $D^b(X)$  and also
$$
\Phi \circ \Phi \cong (-1)^* \id_{D^b(X)}[-3].
$$
The Chern character of any $E \in D^b(X)$ is of the form $\ch(E) = (a_0, a_1 \ell, a_2 \frac{\ell^2}{2}, a_3 \frac{\ell^3}{6})$ for some integers $a_i$.
Then we have (see \cite[Lemma 9.23]{Huy1}):
$$
\ch(\Phi(E)) = (a_3, -a_2 \ell, a_1 \frac{\ell^2}{2}, -a_0 \frac{\ell^3}{6}).
$$

\section{Minimal objects of $\A_{\om , B}$ and B-G Type Inequality of Threefolds}
\label{section2}
\subsection{Some minimal objects of $\A_{\om , B}$}
We identify some classes of minimal objects of the abelian category $\A_{\om,B}$ of a projective threefold $X$.
See \cite{Huy2} for a detailed discussion on minimal objects of some abelian categories associated to Bridgeland stability conditions on a surface.
\begin{prop}
\label{prop2.1}
For any $x \in X$, the skyscraper sheaf $\OO_x$ is a minimal object in $\A_{\om , B}$.
\end{prop}
\begin{proof}
For any $x \in X$, $\OO_x \in \T_{\om , B}$ and also $\OO_x \in \T_{\om , B}'$. Therefore $\OO_x \in \A_{\om , B}$. Let
$$
0 \to a \to \OO_x \to b \to 0
$$
be a short exact sequence (SES for short) in $\A_{\om , B}$ such that $a \ne 0$. Then
in order to prove $\OO_x \in \A_{\om , B}$ is minimal,
we need to show $b =0$. We obtain the following long exact sequence (LES for short) of $\B_{\om , B}$-cohomologies associated to the above $\A_{\om,B}$-SES:
$$
0 \to A_{-1} \to 0 \to B_{-1} \to A_0 \to \OO_x \to B_0 \to 0.
$$
Here $A_k := H^k_{\B_{\om , B}}(a)$ and $B_k := H^k_{\B_{\om , B}}(b)$.
We have $A_{-1} =0$ and so $a \cong A_0 \ne 0$. Let $C: = A_0 / B_{-1}$. Then
$$
0 \to C \to \OO_x \to B_0 \to 0
$$
is a SES in $\B_{\om , B}$. We obtain the following LES of $\coh(X)$-cohomologies associated to the above $\B_{\om,B}$-SES:
$$
0 \to C^{-1} \to 0 \to B_0^{-1} \to C^0 \to \OO_x \to B_0^0 \to 0.
$$
Here $C^k := H^k_{\coh(X)}(C)$ and $B_0^k := H^k_{\coh(X)}(B_0) $. We have  $C^{-1} = 0$ and so $C \cong C^0$.

If $B_0^0 \ne 0$ then $\OO_x \cong B_0^0$ and $ B_0^{-1} \cong C^0 \in
\T_{\om , B} \cap  \F_{\om , B} = \{0\}$. So $C = 0$ and $B_{-1}
\cong A_0 \in  \T_{\om , B}' \cap  \F_{\om , B}' = \{0\}$ which
implies $A_0 = 0$. This is not possible and so $B_0^0 = 0$. Therefore
$B_0 \cong B_0^{-1}[1]$ and
$$
0 \to B_0^{-1} \to C^0 \to \OO_x \to  0
$$
is a SES in $\coh(X)$. Here $\ch(\OO_x) =(0,0,0,1)$. If $B_0^{-1} \ne 0$ then
$$
0 \ge \mu_{\om , B}(B_0^{-1}) =  \mu_{\om , B} (C^0) > 0.
$$
This is not possible and so  $B_0^{-1} = 0$ and  $C^0 \cong \OO_x$. Therefore $b \cong B_{-1}[1]$ and we have the following SES in  $\B_{\om , B}$:
$$
0 \to B_{-1} \to A_0 \to \OO_x \to 0.
$$
Since $\ch(\OO_x) =(0,0,0,1)$, if $B_{-1} \ne 0$ then
$$
0 \ge \nu_{\om , B}(B_{-1}) =  \nu_{\om , B}(A_0) > 0.
$$
This is not possible and so $B_{-1} = 0$. Therefore $b = 0$ and so $\OO_x \in \A_{\om , B}$ is a minimal object as required.
\end{proof}
We now identify further minimal objects.
\begin{defn}
\label{note2.0}
Let $\SM_{\om, B}$ be the class of all objects $E \in \B_{\om , B}$ such that
\begin{enumerate}[label=(\roman*)]
\item $E$ is $\nu_{\om , B}$-stable,
\item $\nu_{\om , B}(E) = 0$, and
\item $\Ext^1(\OO_x , E) = 0$ for any skyscraper sheaf $\OO_x$ of $x \in X$.
\end{enumerate}
\end{defn}
Then clearly $\SM_{\om, B} \subset \SC_{\om,B}$.
\begin{lem}
\label{prop2.2}
Let $E \in \SM_{\om , B}$. Then $E[1]$ is a minimal object of $\A_{\om , B}$.
\end{lem}
\begin{proof}
By definition  $\SM_{\om , B} \subset \F_{\om , B}'$ and so $E[1] \in \A_{\om , B}$.
Let
$$
0 \to a \to E[1] \to b \to 0
$$
be a SES in  $\A_{\om , B}$ such that $b \ne 0$. Now we need to show that $a=0$ or equivalently $b \cong E[1]$.
We have the following LES of $\B_{\om , B}$-cohomologies associated to the above $\A_{\om , B}$-SES:
$$
0 \to A_{-1} \to E \to B_{-1} \to A_0 \to 0 \to B_0 \to 0.
$$
Here $A_k := H^k_{\B_{\om , B}}(a)$ and $B_k := H^k_{\B_{\om , B}}(b)$.
We have $B_0 = 0$ and so $b \cong B_{-1}[1]$ which implies $B_{-1} \ne 0$.
\begin{enumerate}[label={Case} (\Roman*),align=left,labelwidth=*]

\item $A_{-1} \ne 0$:

\begin{enumerate}[label={Subcase} (\roman*),align=left,labelwidth=*]

\item $E/ A_{-1} \ne 0$:\\
Then $E/ A_{-1} \hookrightarrow B_{-1}$ and $\nu_{\om , B}^+(B_{-1}) \le 0$ implies $\nu_{\om , B}(E/ A_{-1}) \le 0$.
On the other hand $\nu_{\om , B}(E/ A_{-1}) > 0$ as $A_{-1} \ne 0$ and $E$ is $\nu_{\om, B}$-stable with $\nu_{\om, B}(E)=0$.
But this not possible.

\item $E/ A_{-1}  =0$:

Then $A_{-1} \cong E$ and $B_{-1} \cong A_0 \in \F_{\om , B}' \cap \T_{\om , B}' = \{0\}$. This is not possible as $B_{-1} \ne 0$.
\end{enumerate}

\item $A_{-1} = 0$:

Then we have the following SES in  $\B_{\om , B}$:
\begin{equation}
0 \to E \to B_{-1} \to A_0 \to 0. \label{ses2.1} \tag{\ding{70}}
\end{equation}
\begin{enumerate}[label={Subcase} (\roman*),align=left,labelwidth=*]
\item $A_{0} \ne 0$:

Here $\nu_{\om, B}(E) = 0$ implies $\om^2 \ch_1^B(E) > 0$ and $\Im\,Z_{\om, B}(E) =0$. Then
$$
\nu_{\om, B}(B_{-1}) = \frac{\Im\,Z_{\om, B}(A_0)}{\om^2 \ch_1^B(E) + \om^2 \ch_1^B(A_0)} \le 0
$$
implies $\Im\,Z_{\om, B}(A_0) \le 0$.
If $\om^2 \ch_1^B(A_0) \ne 0$ then $\nu_{\om, B}(A_0) > 0$ implies $\Im\,Z_{\om, B}(A_0) >0 $; which is not possible. Hence $\om^2 \ch_1^B(A_0)=0$ and by Lemma \ref{prop1.1}, $\Im\,Z_{\om, B}(A_0) \ge 0$. So $ \Im\,Z_{\om, B}(A_0) = 0$ and $A_0 \cong T$ for  some $0 \ne T \in \coh^0(X)$.
Then the $\B_{\om, B}$-SES \eqref{ses2.1}  corresponds to an element from $\Ext^1(A_0 , E) = \Ext^1(T , E)$. But we have $\Ext^1(\OO_x , E) = 0$ for any $x \in X$ and so  $\Ext^1(T , E) =0$. So $ B_{-1} \cong T \oplus E $. Then
$T$ is a subobject of $B_{-1}$. But this is not possible as $\nu_{\om, B}(T) = + \infty$ and $E \in \SM_{\om,B}$.

\item $A_{0} = 0$:

Then $a = 0$ and $ b \cong B_{-1}[1] \cong E[1]$ as required.
\end{enumerate}
\end{enumerate}
This completes the proof of the lemma.
\end{proof}
Some classes of tilt stable candidates have been identified in \cite{BMT}.

Recall, for $E \in D^b(X)$ the discriminant $\overline{\Delta}_\om$ in the sense of Dr\'{e}zet is defined by
$$
\overline{\Delta}_\om(E)  = \left( \om^2 \ch_1^B(E) \right)^2 -  2 \om^3 \ch_0^B(E) \cdot \om \ch_2^B(E).
$$
\begin{prop} \cite[Proposition 7.4.1]{BMT}
\label{prop2.3}
Let $E$ be a $\mu_{\om, B}$-stable locally free sheaf on $X$ with $\overline{\Delta}_\om(E) =0$. Then either $E$ or $E[1]$ in $\B_{\om, B}$ is  $\nu_{\om, B}$-stable.
\end{prop}
\begin{exmp}
\label{ex2.1}
Let $(X,L)$ be a polarized projective threefold and let $\ell:= c_1(L)$. Consider the classes $B = \frac{1}{2} \ell$ and $\om = \frac{\sqrt{3}}{2} \ell$.
Then  $\overline{\Delta}_\om(\OO) = \overline{\Delta}_\om(L) = 0$. So, by Proposition \ref{prop2.3}, $\OO[1] , L \in \B_{\om, B}$ are $\nu_{\om, B}$-stable. Also $\Im\,Z_{\om, B}(\OO[1]) = \Im\,Z_{\om, B}(L) = 0$.  Therefore $\nu_{\om, B}(\OO[1]) = \nu_{\om, B}(L) = 0$.
So by Lemma \ref{prop2.2}, $\OO[2] , L[1] \in \A_{\om, B}$ are minimal objects.
\end{exmp}
\begin{note}
\label{note2.1}
The tilt stable objects associated to minimal objects in Example
\ref{ex2.1} clearly satisfy the corresponding B-G type inequalities.
\end{note}
\subsection{Reduction of B-G type inequality for minimal objects}

The following propositions are important.
\begin{prop} \cite[Proposition 3.1]{LM}
\label{prop2.4}
Let $E \in \B_{\om ,B}$ be a $\nu_{\om,B}$-semistable object with $\nu_{\om, B} < +\infty$. Then
$H^{-1}_{\coh(X)}(E)$ is a reflexive sheaf.
\end{prop}
\begin{prop} \cite[Proposition 3.5]{LM}
\label{prop2.5}
Let
$0 \to E \to E' \to Q \to 0$
be  a non splitting SES in $\B_{\om,B}$ with
 $Q \in \coh^{0}(X)$, $\Hom(\OO_x, E') = 0$ for all $x \in X$, and
 $\om^2 \ch_1^B(E) \ne 0$. If $E$ is $\nu_{\om, B}$-stable then $E'$
 is $\nu_{\om, B}$-stable.
\end{prop}
Recall that $\SC_{\om, B}$ is the class of $\nu_{\om , B}$-stable
objects  $E \in \B_{\om , B}$ with $\nu_{\om , B}(E) = 0$.
\begin{prop}
\label{prop2.6}
Let $E \in  \SC_{\om, B}$. Then there exists $E' \in \SM_{\om, B}$
(i.e. $E'[1]$ is a minimal object of $\A_{\om, B}$)  such that
$$
0 \to E \to E' \to Q \to 0
$$
is a SES in $\B_{\om, B}$ for some $Q \in \coh^{0}(X)$.
\end{prop}
\begin{proof}
Let $E \in  \SC_{\om, B} \setminus \SM_{\om, B}$.
Assume the opposite of the claim in the proposition for $E$. Then
there exists a sequence of non-splitting SESs in $\B_{\om,B}$, for $i
\ge 1$
$$
0 \to E_{i-1} \to E_i \to \OO_{y_i} \to 0,
$$
where $E_0 = E$, $E_i \in \SC_{\om, B}$ (see Proposition \ref{prop2.5}).
So for each $i \ge 1$,
$$
0 \to \OO_{y_i} \to E_{i-1}[1] \to E_i[1] \to 0
$$
is a SES in $\A_{\om, B}$. Therefore
$$
E[1] = E_0[1] \twoheadrightarrow E_1[1] \twoheadrightarrow E_2[1] \twoheadrightarrow \cdots
$$
is an infinite chain of quotients in $\A_{\om,B}$. But this is not
possible as $\A_{\om, B}$ is Noetherian by Lemma \ref{prop1.4}. This
is a contradiction.
\end{proof}
It follows that $E \in \SC_{\om,B}$ satisfies the B-G type inequality
if the corresponding $E' \in \SM_{\om,B}$ satisfies the B-G type
inequality.
\section{Abelian category $\A_{\sqrt{3}B, B}$, FMT and stability conditions}
\label{section3}
\subsection{Some properties of $\A_{\sqrt{3}B, B}$}
We discuss some of the properties of the abelian category
$\A_{\sqrt{3}B, B}$ for an arbitrary polarized projective threefold
$(X,L)$ with Picard rank $1$. Let $\ell:= c_1(L)$. Let $B = b \ell$
for $b \in \Q_{>0}$. Then for $E \in D^b(X)$
\[
\Im\,Z_{\sqrt{3}B, B}(E) =  \sqrt{3} b \ell ( \ch_2(E) - b \ell \ch_1(E)).
\]
\begin{prop}
\label{prop3.1}
Let $E \in \B_{\sqrt{3}B, B}$ and  let $E_{i} =
H^{i}_{\coh(X)}(E)$. Let $E_i^\pm$ be the  H-N semistable factors of
$E_i$ with highest and lowest $\mu_{\sqrt{3}B,B}$ slopes. Then we have
the following:
\begin{enumerate}[label=(\roman*)]
\item if $E \in \HN^{\nu}_{\sqrt{3}B , B}(-\infty, 0) $ and $E_{-1} \ne 0$, then $\ell^2 \ch_1(E_{-1}^+) < 0$;
\item if $E \in \HN^{\nu}_{\sqrt{3}B , B}(0, +\infty]$ and $\rk(E_0) \ne 0$, then $\ell^2 \ch_1(E_0^-) > 2 b \ell^3 \ch_0(E_0^-)$; and
\item if $E$ is tilt semistable with $\nu_{\sqrt{3}B, B}(E) =0$, then
\begin{enumerate}
\item for $E_{-1} \ne 0$, $\ell^2 \ch_1(E_{-1}) \le  0$ with equality
  if and only if $\ch_2(E_{-1}) = 0$,
\item for $\rk(E_0) \ne 0$, $\ell^2 \ch_1(E_0) \ge 2 b \ell^3
  \ch_0(E_0)$ with equality
  if and only if $(\ch_1(E_0))^2 = 2 \ch_0(E_0) \ch_2(E_0)$.
\end{enumerate}
\end{enumerate}
\end{prop}
\begin{proof}
$E \in  \B_{\sqrt{3}B, B}$ fits in to the $\B_{\sqrt{3}B, B}$-SES
$$
0 \to E_{-1}[1] \to E \to E_0 \to 0.
$$
\begin{enumerate}[label=(\roman*)]
\item
Since $E \in \HN^{\nu}_{\sqrt{3}B , B}(-\infty, 0)$, $E_{-1}[1] \in \HN^{\nu}_{\sqrt{3}B , B}(-\infty, 0)$. We have $0 \ne E_{-1}^+ \subseteq E_{-1}$.
Hence $E_{-1}^+[1] \in \HN^{\nu}_{\sqrt{3}B , B}(-\infty, 0)$.

Let $\ch(E_{-1}^+)=(a_0, a_1, a_2, a_3)$. Assume the opposite for a contradiction; so that
$a_1 \ge 0$.
We have
\begin{align*}
  \nu_{\sqrt{3}B, B} (E_{-1}^+[1])   & = \frac{-\Im\,Z_{\sqrt{3}B, B}(E_{-1}^+)}{-3B^2 \ch_1^B(E_{-1}^+)} \\
                               & = \frac{\sqrt{3}b a_1(b a_0-a_1) + \frac{\sqrt{3}}{2}b a_1^2 + \frac{\sqrt{3}}{2}b ( a_1^2 - a_0 a_2)}{3a_0b^2(ba_0-a_1)}.
\end{align*}
Since $E_{-1}^+$ is $\mu_{\sqrt{3}B, B}$-semistable we have, by the usual B-G inequality,
$$
a_1^2 -a_0 a_2  \ge 0
$$
and since $E_{-1}^+\in\F_{\sqrt{3}B,B}$, we have $\nu_{\sqrt{3}B, B} (E_{-1}^+[1]) \ne +\infty$ and so $ba_0-a_1 > 0$.
Hence, as $a_0>0$, we have $\nu_{\sqrt{3}B, B}
(E_{-1}^+[1]) \ge 0$. But this is not possible as $E_{-1}^+[1] \in
\HN^{\nu}_{\sqrt{3}B , B}(-\infty, 0)$. This is the required
contradiction to complete the proof.

\item
Since $E \in \HN^{\nu}_{\sqrt{3}B , B}(0, +\infty]$, $E_0 \in \HN^{\nu}_{\sqrt{3}B , B}(0, +\infty]$.
We have $0 \ne E_{0}^-$ is a torsion free quotient of $E_0$. Since
$E_0 \in \HN^{\nu}_{\sqrt{3}B , B}(0, +\infty]$
we have $E_0^- \in \HN^{\nu}_{\sqrt{3}B , B}(0, +\infty]$ .

Let $\ch(E_{0}^-)=(a_0, a_1, a_2, a_3)$. Assume the opposite for a contradiction; so that $a_1 \le 2b a_0$.
We have
\begin{align*}
\nu_{\sqrt{3}B, B} (E_0^-)  & =   \frac{\Im\,Z_{\sqrt{3}B, B}(E_0^-)}{3B^2 \ch_1^B(E_0^-)} \\
                      & = \frac{ -\frac{\sqrt{3}}{2}b (a_1^2 -a_0 a_2) +  \frac{\sqrt{3}}{2}b a_1 (a_1 - 2ba_0)}{3b^2 a_0(a_1-ba_0)}.
\end{align*}
Here $E_0^- \in \T_{\sqrt{3}B, B}$ is torsion free which implies
$$
 a_1 - b a_0  > 0;
$$
$E_0^-$ is $\mu_{\sqrt{3}B, B}$-semistable which implies (by the usual B-G inequality)
$$
a_1^2 - a_0 a_2 \ge 0.
$$
Therefore $\nu_{\sqrt{3}B, B} (E_0^-) \le 0$. But this is not possible as $E_0^- \in \HN^{\nu}_{\sqrt{3}B , B}(0, +\infty]$.
This is the required contradiction to complete the proof.

\item
Similar to (i) one can show that if $E \in \HN^{\nu}_{\sqrt{3}B , B}(-\infty, 0] $ and $E_{-1} \ne 0$, then $\ell^2 \ch_1(E_{-1}^+) \le 0$.
Therefore for $E \in \HN^{\nu}_{\sqrt{3}B , B}[0]$ we have  $\ell^2 \ch_1(E_{-1}) \le  0$.
The equality holds if and only if $E_{-1}$ is slope semistable, and so it satisfies the usual B-G inequality.
Since $\nu_{\sqrt{3}B, B}(E_{-1}) \le 0$  we have $\ell^2 \ch_1(E_{-1}) =  0$  if and only if $\ch_2(E_{-1}) = 0$.

Proof of (b) is similar to that of (a).
 \end{enumerate}
\end{proof}
\subsection{Relation of FMT to stability conditions}
Let $(X, L)$ be a principally polarized abelian threefold with Picard rank 1.
Let $\ell: = c_1(L)$. Then $\chi(L) = \frac{\ell^3}{6}= 1$ and the Chern character of  $E \in D^b(X)$ is of the form $\ch(E) = (a_0, a_1 \ell, a_2 \frac{\ell^2}{2}, a_3 \frac{\ell^3}{6})$ for some integers $a_i$.
Define the classes  $B  = \frac{1}{2} \ell$ and $\omega = \frac{\sqrt{3}}{2} \ell$.

The following is a key result in this paper.
\begin{prop}
\label{prop3.2}
If $\Phi (L^{-1} E)[2] \in \B_{\om, B}$ for any $E \in \SM_{\om, B} \setminus \{ L \Po_x : x \in X \}$, then
the B-G type inequality holds for the objects in $\SC_{\om,B}$.
\end{prop}
\begin{proof}
By Proposition \ref{prop2.6}, it is enough to check that the B-G type inequality is satisfied by each object in $\SM_{\om,B}$. Moreover,  the objects  in $\{ L \Po_x : x \in X \} \subset \SM_{\om,B}$ satisfy the B-G type inequality (see Note \ref{note2.1}). Then we only need to check the inequality for objects in $\SM_{\om, B} \setminus \{ L \Po_x : x \in X \}$.

Let  $E \in \SM_{\om, B} \setminus \{ L \Po_x : x \in X \}$ and assume $\Phi (L^{-1} E)[2] \in \B_{\om, B}$.
Let $\ch(E) = (a_0, a_1 \ell, a_2 \frac{\ell^2}{2}, a_3\frac{\ell^3}{6})$ and then
$\Im\,Z_{\om,B}(E) = 0$ implies $a_1 = a_2$. Now the B-G type inequality says
\begin{equation*}
\Delta : = -a_0 + 3a_1 - a_3 > 0.
\end{equation*}

By Proposition \ref{prop3.1}, we have $\ell^2 \ch_1(E_{-1}) \le 0$ and
$\ell^2 \ch_1(E_0) \ge 0$. Here $E_i = H^{i}_{\coh(X)}(E)$. So $a_1 \ell^3 = \ell^2 \ch_1(E) = \ell^2 \ch_1(E_0) -  \ell^2 \ch_1(E_{-1}) \ge 0$.

Let $F  = \Phi (L^{-1} E)[2]$ and let $\ch(F) = (b_0, b_1 \ell, b_2 \frac{\ell^2}{2}, b_3 \frac{\ell^3}{6})$. Then
$b_0 = a_3 - a_0$ and $b_1 = b_2 = a_1 -a_0$. Now $b_1 = b_2$ implies $\Im\,Z_{\om ,B}(F) = 0$. Also $F \in  \B_{\om, B}$ implies
$\om^2 \ch_1^B(F) \ge 0$, i.e. $2b_1 - b_0 \ge 0$.
If  $\om^2 \ch_1^B(F) = 0$  then $\Im\,Z_{\om ,B}(F) =0$ implies $F \cong T$ for some $T \in \coh^0(X)$ (see Lemma \ref{prop1.1}).
If $T\neq 0$ then $E$ has a filtration with factors of the form $L
\Po_x[1]\not\in\SM_{\om,B}$. This is not possible and so $\om^2
\ch_1^B(F) > 0$. That is $2b_1 - b_0  = -a_0 + 2a_1 -a_3 > 0$.

Hence $\Delta > 0$ and so $E$ satisfies the B-G type inequality. This completes the proof as required.
\end{proof}
Our main goal in the rest of this paper is to prove that $\Phi
L^{-1}[2]$ and its quasi-inverse $L\Phi[1]$ are auto-equivalences of
the abelian category $\A_{\om,B}$.
Under an equivalence of abelian categories minimal objects are mapped
to minimal objects and so the hypothesis of Proposition
\ref{prop3.2} is satisfied. Therefore, by Corollary \ref{prop1.5}, we
have the following:

\begin{thm}
\label{prop3.3}
The pair $(\A_{\om, B}, Z_{\om, B})$ is a Bridgeland stability condition on $D^b(X)$.
\end{thm}

\section{Fourier-Mukai transforms on $\mathrm{Coh}(X)$ of Abelian Threefolds}
\label{section4}

From here onward, we always assume $(X, L)$ is a principally polarized abelian threefold with Picard rank 1.
Let $\ell: = c_1(L)$. Then $\chi(L) = \frac{\ell^3}{6}= 1$ and the Chern character of any $E \in D^b(X)$ is of the form $\ch(E) = (a_0, a_1 \ell, a_2 \frac{\ell^2}{2}, a_3 \frac{\ell^3}{6})$ for some integers $a_i$. Define the classes  $B  = \frac{1}{2} \ell$ and $\omega = \frac{\sqrt{3}}{2} \ell$.

If $E \in \coh(X)$ then the slope $\mu(E)$ is defined by $\mu(E) := \mu_{\frac{1}{\sqrt{6}}\ell, 0}(E)$.
That is $\mu(E) = \frac{a_1}{a_0}$ when $a_0 \ne 0$, and $\mu(E) = +\infty$ when $a_0 = 0$. In the rest of the paper we mostly use $\mu$ slope for coherent sheaves and we simply write
$\HN = \HN^{\mu}_{\frac{1}{\sqrt{6}}\ell, 0}$.
Then  $\mu_{\om, B} (E) = \frac{9}{2}(\mu(E) - \frac{1}{2})$.
Moreover define $\T_0 = \HN(0,+\infty]$ and $\F_0 = \HN(-\infty, 0]$.
Also for simplicity we write
$\T = \T_{\om, B}$, $\F = \F_{\om, B}$, $\B = \B_{\om, B}$, $\nu = \nu_{\om, B}$, $\HN^{\nu}_{\om, B} = \HN^{\nu}$,
$\T' = \T_{\om, B}'$, $\F' = \F_{\om, B}'$, and $\A = \A_{\om, B}$.
Then by the definitions, we have $\F =\HN(-\infty, \frac{1}{2}]$ and $\T = \HN(\frac{1}{2}, + \infty]$.

Let $\Phi$ be the Fourier-Mukai transform with kernel the Poincar\'{e}
line bundle $\Po$. The isomorphism $\Phi \circ \Phi \cong (-1)^*
\id_{\id_{D^b(X)}}[-3]$ gives us the following convergence of spectral
sequence.
\begin{mukaispecseq}
\label{mukaiss}
$$
E_2^{p,q} = \Phi^p_{\coh(X)}  \Phi^q_{\coh(X)}(E) \Longrightarrow H^{p+q-3}_{\coh(X)}((-1)^* E),
$$
for $E$. Here $\Phi^i_{\coh(X)} (F) = H^i_{\coh(X)}(\Phi(F))$.
\end{mukaispecseq}
 For $E \in \coh(X)$, we write
$$
E^k = \Phi^k_{\coh(X)}(E).
$$
Then for example $E^{120} =\Phi^0_{\coh(X)} \Phi^2_{\coh(X)}\Phi^1_{\coh(X)}(E)$.
Using this notation, we can deduce the following immediately from the spectral sequence:
\[E^{00}=E^{01}=E^{32}=E^{33}=0\text{, }E^{10}\cong E^{02}\text{ and
}E^{31}\cong E^{23}.\]
Let $\RR \Delta$ denote the derived dualizing functor $\RR \calHom(-, \OO)[3]$. Then due to Mukai,
\[(\Phi \circ \RR \Delta) [3]\cong(-1)^*\RR \Delta \circ \Phi\]
(see \cite[(3.8)]{Muk2}).
This gives us the convergence of the following spectral sequences.
\begin{dualspecseq}
\label{dualss}
$$
\Phi^p_{\coh(X)} \left( \calExt^{q+3} (E, \OO) \right) \Longrightarrow \ \ ? \ \ \Longleftarrow (-1)^* \calExt^{p+3}\left( \Phi^{3-q}_{\coh(X)}(E), \OO \right)
$$
for $E \in \coh(X)$.
\end{dualspecseq}
The aim of this section is to use mainly the Mukai and ``Duality'' spectral sequences 
to study the slope stability of sheaves under the FM transform $\Phi$.
More precisely, we consider the $\coh(X)$ cohomology sheaves of the
images under $\Phi$ of torsion sheaves supported in 
dimensions $1$ and $2$. We also study the transforms of torsion free
sheaves whose H-N semistable factors satisfy certain slope bounds. 
\begin{notation}
Any $E \in \coh(X)$ fits into $\coh(X)$-SES
$$
0 \to T \to E \to F \to 0
$$
for some $T \in \T_0$ and $F \in \F_0$. Denote $T(E)  = T$ and $F(E) = F$.
\end{notation}
Any torsion free sheaf $E$ fits into a non-splitting $\coh(X)$-SES
$$
0 \to E \to E^{**} \to T \to 0
$$
for some $T \in \coh^{\le 1}(X)$. Here $E^{**}$ is a reflexive sheaf.
If $E$ is rank 1 then $E^{**}$ is a line bundle and so $E^{**} \cong L^k \Po_x$ for some $k \in \Z$ and $x \in X$.
\begin{notation}
\label{prop4.1}
If $E$ is a rank 1 torsion free sheaf with $c_1 (E) = k \ell$ then we
can write $E = L^k \Po_x \I_C$. Here $\I_C$ is the ideal sheaf of the
structure sheaf $\OO_C : = L^{-k}\Po_{-x}\otimes(E^{**} / E) \in \coh^{\le
  1}(X)$ of a subscheme $C \subset X$  of dimension $\le 1$.
\end{notation}
\begin{prop}
\label{prop4.2}
Let $E \in \coh(X)$. If $E^0 \ne 0$ then $E^0$ is a reflexive sheaf.
\end{prop}
\begin{proof}
Let $x \in X$. Then for $0 \le i \le 2$, we have
$$
\Hom(\OO_x, E^0[i]) \cong \Hom(\Phi(\OO_x), \Phi(E^0)[i]) \cong \Hom(\Po_x, E^{02}[-2+i])
$$
from the convergence of the Mukai Spectral Sequence \ref{mukaiss} for $E$.
So $\Hom(\OO_x, E^0) = \Ext^1(\OO_x, E^0) =  0$, and
\begin{align*}
\Ext^2(\OO_x, E^0) & \cong \Hom (\Po_x, E^{02}) \\
                   & \cong \Hom(\Po_x, E^{10}),  \ \ \mbox{by the Mukai Spectral Sequence for $E$}\\
                   & \cong \Hom(\Phi(\OO_x), \Phi(E^1)) \\
                   & \cong \Hom(\OO_x, E^1).
\end{align*}
Hence, as any map $\OO_x\to E^1$ must factor through the torsion
subsheaf of $E^1$ and $E^1$ is coherent, only finitely many of these
can be non-zero. So $\dim \{ x \in X:
\Ext^2(\OO_x, E^0 ) \ne 0 \} \le 0$. Therefore $E^0$ is a reflexive
sheaf.
\end{proof}
\begin{prop}
\label{prop4.3}
Let $E \in \coh(X)$. Then we have the following:
\begin{enumerate}[label=(\roman*)]
\item if $E \in \T_0$ then $E^3 = 0$, and
\item if $E \in \F_0$ then $E^0 = 0$.
\end{enumerate}
\end{prop}
\begin{proof}
\begin{enumerate}[label=(\roman*)]
\item
Let $E \in \T_0$. Then for any $x \in X$, we have
\begin{align*}
\Hom(E^3, \OO_x) & \cong \Hom(\Phi(E)[3], \Phi(\Po_{-x})[3]) \\
                 & \cong \Hom(E, \Po_{-x}) = 0,
\end{align*}
as $\Po_{-x} \in \F_0$. Therefore $E^3 =0$ as required.
\item
Let $E \in \F_0$. We can assume $E$ is $\mu$-stable using H-N and
Jordan-H\"older filtrations. For generic $x \in X$ and $i=1,2$ we have
$$
\Hom(E^1, \OO_x[i]) = \Hom(E^2, \OO_x[i+1]) = \Hom(E^3, \OO_x[i+2]) = 0.
$$
Hence for generic $x \in X$,
\begin{align*}
\Hom(E^0, \OO_x) & \cong \Hom(\Phi(E), \OO_x) \\
                 & \cong \Hom(\Phi(E), \Phi(\Po_{-x})[3]) \\
                 & \cong \Hom(E, \Po_{-x}[3]) \\
                 & \cong \Hom (\Po_{-x}, E)^*.
\end{align*}
\begin{enumerate}[label=(\alph*)]
\item Case $\mu(E) < 0$: \\
Then $\Hom(\Po_{-x}, E) = 0$.
\item Case $\mu(E) =0$: \\
Since $E$ is assumed to be $\mu$-stable, any map in
$\Hom(\Po_{-x},E)$ must be an isomorphism and so $E^0=0$.
\end{enumerate}
Therefore for generic $x \in X$, $\Hom(E^0, \OO_x)= 0$. By Proposition
\ref{prop4.2} if  $E^0 \ne 0$ then it is reflexive. So $E^0 = 0$.
\end{enumerate}
\end{proof}
\begin{prop}
\label{prop4.4}
Let $E \in \coh(X)$. Then
\begin{enumerate}[label=(\roman*)]
\item  $E^3 \in \T_0$, and
\item  $E^0 \in \F_0$.
\end{enumerate}
\end{prop}
\begin{proof}
\begin{enumerate}[label=(\roman*)]
\item
Let $T = T(E^3) \in \T_0$ and $F = F(E^3) \in \F_0$,  so that
$$
0 \to T \to E^3\to F \to 0
$$
is a SES in $\coh(X)$. Now we need to show that $F=0$.
Apply $\Phi$ to the above SES and consider the LES of $\coh(X)$-cohomologies. Then we have
$F \in V^{\Phi}_{\coh(X)}(1)$, $T \in V^{\Phi}_{\coh(X)}(0,1,2)$ (for
the definition of $V$ see the notation section of the introduction) and
$$
0 \to T^1 \to E^{31} \to F^1 \to T^2 \to 0
$$
is a LES in $\coh(X)$. Here $E^{31} \cong E^{23}$ (from the Mukai Spectral Sequence \ref{mukaiss} for $E$) and so
\begin{align*}
\Hom(E^{31}, F^1) & \cong  \Hom (E^{23}, F^1) \\
                  & \cong  \Hom (\Phi (E^2) [3] , \Phi (F) [1]) \\
                  & \cong  \Hom (E^2 , F[-2]) = 0.
\end{align*}
Hence $F \cong (-1)^*F^{12} \cong  (-1)^* T^{22} = 0$ (from the Mukai Spectral Sequence \ref{mukaiss} for $T$) as required.

\item Similar to the proof of (i).
\end{enumerate}
\end{proof}
\begin{prop}
\label{prop4.6}
Let $E \in \F_0$. If $E^1 \ne 0$ then $E^1$ is a reflexive sheaf.
\end{prop}
\begin{proof}
By Proposition \ref{prop4.3}, $E^0 = 0$.
Let $x \in X$. Then from the convergence of the Mukai Spectral Sequence \ref{mukaiss} for $E$ and $0 \le i \le 2$, we have
\begin{align*}
 \Hom(\OO_x, E^1[i]) & \cong  \Hom(\Phi(\OO_x), \Phi(E^1)[i]) \\
                     & \cong  \Hom(\Po_x, E^{12}[i-2])
\end{align*}
as $\Hom(\Po_x,\tau_{>2}\Phi(E^1)[i])\cong\Hom(\Po_x,E^{13}[i-3])=0$. Therefore $\Hom(\OO_x, E^1) =\Ext^1(\OO_x, E^1)=0$ and $\Ext^2(\OO_x, E^1)  \cong \Hom(\Po_x , E^{12})$.

From the convergence of the Mukai Spectral Sequence \ref{mukaiss} for $E$
$$
0 \to E^{20} \to E^{12} \to F \to 0
$$
is a SES in $\coh(X)$. Here $F$ is a subobject of $(-1)^*E$. By
applying the functor $\Hom(\Po_x, -)$ we obtain the  exact sequence
$$
0 \to \Hom(\Po_x, E^{20}) \to \Hom(\Po_x, E^{12}) \to \Hom(\Po_x , F) \to \cdots.
$$
Now $F\in \F_0$ and by Proposition \ref{prop4.4} $E^{20}$ is also in $\F_0$. Therefore
we have $\Hom(\Po_x , F) \ne 0$ or $\Hom(\Po_x, E^{20}) \ne 0$ for at
most a finite number of points $x \in X$. That is
$\dim \{ x \in X: \Ext^2(\OO_x, E^1) \ne 0 \} \le 0$. Therefore $E^1$ is a reflexive sheaf.
\end{proof}
\begin{prop}
\label{prop4.7}
If $E$ is a torsion sheaf then $E^2 \in \T_0$.
\end{prop}
\begin{proof}
Let $T = T(E^2)$ and $F= F(E^2)$. Then $0 \to T \to E^2 \to F \to 0$ is a SES in $\coh(X)$. By applying $\Phi$ we obtain the  LES
$$
0 \to T^1 \to E^{21} \to F^1 \to T^2 \to 0
$$
in $\coh(X)$. Here $F \in V^{\Phi}_{\coh(X)}(1)$.
From the convergence of the Mukai Spectral Sequence \ref{mukaiss} for $E$, $E^{21}$ fits into the $\coh(X)$-SES
$$
0 \to Q \to E^{21} \to E^{13} \to 0,
$$
where $Q$ is a quotient of $(-1)^*E$.
So $Q$ is a torsion sheaf and $\Hom (Q, F^1) = 0$ as $F^1$ is a reflexive sheaf (see Proposition \ref{prop4.6}). Therefore
\begin{align*}
\Hom (E^{21}, F^1) & \cong  \Hom (E^{13} , F^1) \\
                   & \cong  \Hom (\Phi (E^1) [3] , \Phi (F) [1]) \\
                  & \cong  \Hom (E^1 , F[-2]) = 0.
\end{align*}
Hence $F^1 \cong T^2$ and so $F \cong(-1)^* F^{12} \cong (-1)^*T^{22} = 0$ (from the Mukai Spectral Sequence \ref{mukaiss} for $T$) as required.
\end{proof}
For $x \in X$, let $L_x$ denote $L \Po_x$. 
Since $h^0(X, L_x) = \chi(L_x) = 1$, let the divisor $D_x$ 
be the zero locus of the unique (up to scale) section $0 \ne s_x \in H^0(X, L_x)$.
Moreover, as $t_x^*L \otimes L^{-1} = \Po_x$, we have $D_x = t_x^* D_e$, where $e \in X$ is the identity element.
For positive integer $m$, let $mD_x$ be the non-reduced divisor in the linear system $|m \ell|$
topologically supported on $D_x$.
So $mD_x$ is the zero locus of the section $s_x^{\otimes m}$ of $L_x^m$, and
 we have the SES
\begin{equation*}
\label{ses_res}
0 \to L_x^{-m} {\to} \OO_X \to \OO_{mD_x} \to 0
\end{equation*}
in $\coh(X)$.
For $E \in \coh(X)$, 
apply the functor $E \stackrel{\textbf{L}}{\otimes} (-)$ to  the above SES and consider the LES
of $\coh(X)$-cohomologies.
Since $L_x^{-m}$ and $\OO_X$ are locally free, we have $\Tor_i(E, \OO_{mD_x}) = 0$ for $i \ge 2$.
Now assume  $E \in \coh^k(X)$ for some $k \in \{0,1,2,3\}$.
 For generic $x \in X$, we have $\dim (\Supp(E) \cap {D_x}) \le (k-1)$ and so
 $\Tor_1(E, \OO_{mD_x}) \in \coh^{\le k-1}(X)$.
However, $ L_x^{-m} E \in \coh^k(X)$, and so  $\Tor_1(E, \OO_{mD_x}) =0$.
Therefore, we have the SES
\begin{equation}
\label{ses_res_any}
0 \to L_x^{-m} E \to E \to E|_{mD_x} \to 0
\tag{\ding{61}}
\end{equation}
in $\coh(X)$.
Since any $E \in \coh(X)$ is an extension of sheaves from $\coh^k(X)$,
for generic $x \in X$, $\Tor_i(E, \OO_{mD_x}) = 0$ for $i \ge 1$
and so we have the SES \eqref{ses_res_any}.
Moreover, when $0 \to E_1 \to E_2 \to E_3 \to 0$ is a SES in $\coh(X)$,
for generic $x \in X$ we have $\Tor_i(E_j, \OO_{mD_x}) = 0$, $i \ge 1$ for each $j$, and so
$$
0 \to E_1|_{mD_x} \to E_2|_{mD_x} \to E_3|_{mD_x} \to 0
$$
is a SES in $\coh(X)$.
\begin{prop}
\label{prop4.8}
Let $E \in \coh^{\le 1}(X)$. Then $E^1 \in \T_0$.
\end{prop}
\begin{proof}
$E \in \coh^{\le 1}(X)$ fits into the torsion sequence $0 \to E_0 \to
  E \to  E_1 \to 0$, where $E_0 \in \coh^0(X)$ and $E_1 \in
  \coh^1(X)$. Here $E_0 \in V^{\Phi}_{\coh(X)}(0)$ and so $E^1 =
  E_1^1$. Therefore we only need to prove the claim for a pure
  dimension 1 torsion sheaf $E$. Then  for sufficiently large $m>0$ and suitable $x\in X$, $L_x^{-m} E \in
  V^{\Phi}_{\coh(X)}(1)$, and
$$
0 \to L_x^{-m}E \to E \to E|_{mD_x} \to 0
$$
is a SES in $\coh(X)$ for  $E|_{mD_x} \in \coh^0(X)$. By applying the FMT $\Phi$, we have
$(L_x^{-m}E)^1 \twoheadrightarrow \ E^1$. Therefore, we only need to
show $(L_x^{-m}E)^1 \in \T_0$. Let us show this by proving the claim
for a pure dimension one torsion sheaf $E \in
V^{\Phi}_{\coh(X)}(1)$. Then $\ch(E) = (0,0,\alpha, \beta)$, where
$\alpha> 0$ and $\beta \leq 0$ since $\beta=-\rk(E^1)$.

Let $T = T(E^1)$ and $F = F(E^1)$. Then $ 0 \to  T \to E^1 \to F \to
0$ is a SES in $\coh(X)$.
Now we need to show $F = 0$. So suppose $F \ne 0$ for a contradiction.
Apply the FMT $\Phi$ and consider the LES of $\coh(X)$-cohomologies. Then we have $T \in V^{\Phi}_{\coh(X)}(2)$,  $F \in V^{\Phi}_{\coh(X)}(1,2)$ and
$$
0 \to F^1 \to T^2 \to E \to F^2 \to 0
$$
is a LES in $\coh(X)$.
\begin{enumerate}[label={Case} (\roman*),align=left,labelwidth=*]
\item The map $T^2 \to E$ is zero:

Then $T \cong (-1)^* T^{21} \cong (-1)^* F^{11} = 0$ from the Mukai
Spectral Sequence \ref{mukaiss} as $F \in V^{\Phi}_{\coh(X)}(1,2)$. So $E = F^2$ and
hence $F\in V^{\Phi}_{\coh(X)}(2)$. Therefore $F \cong (-1)^* E^1$ and so
$\ch(F) = ( - \beta, \alpha, 0 ,0 )$.
Here $\alpha > 0$ and which is not possible as $\mu(F) \le 0$.

\item The map $T^2 \to E$ is non-zero:

Let $K = \im (T^2 \to E)$.
Then $K \in \coh^{1}(X)$ and the $\coh(X)$-SES $0 \to F^1 \to T^2 \to K \to 0$ corresponds to an element from
$\Ext^1(K, F^1)$. Here $F^1$ is a reflexive sheaf and so there exists
a locally free sheaf $U$ and a torsion free sheaf $V$ such that
$0 \to F^1 \to U \to V \to 0$ is a non-splitting SES in $\coh(X)$. By
applying the functor $\Hom(K, -)$, we obtain the following exact
sequence:
$$
\cdots \to \Hom (K, V) \to \Ext^1(K, F^1) \to \Ext^1(K, U) \to \cdots.
$$
Here $\Hom (K, V) =0 $ and $\Ext^1(K, U) \cong \Ext^2(U, K)^*\cong
H^2(X,U^*\otimes K)^* = 0$ as $K\in\coh^{\leq1}(X)$.
So $\Ext^1(K, F^1) = 0$ implies $T^2 \cong F^1 \oplus K$. Here $T^2 \in
V^{\Phi}_{\coh(X)}(1)$ implies $F^1 = 0$ and so $K \cong T^2$. Then $F^2 \cong
E/T^2 $ and also $F \in V^{\Phi}_{\coh(X)}(2)$.
Since $F^2 \in V^{\Phi}_{\coh(X)}(1)$, it is  a pure dimension 1 torsion sheaf. So $\ch(F^2) = (0,0, \alpha', \beta')$, where
$\alpha' > 0$ and $\beta'\leq 0$. Therefore $\ch(F) = (-\beta', \alpha', 0, 0)$ and which is not possible as $\mu(F) \le 0$ implies $\alpha' \le 0$.
\end{enumerate}
Therefore $F = 0$ as required to complete the proof.
\end{proof}
Recall from \cite[Prop 6.16]{Muk1} for any positive
 integer $s$, the semi-homogeneous bundle $(L^s)^0$ is slope stable.
In the rest of this section we abuse notation to write $(L^s)^0$ for the
functor $(L^s)^0 \otimes -$.
\begin{prop}
\label{prop4.9.1}
Let $E_n \in \HN[0, +\infty)$, $n \in \N$ be a sequence of coherent sheaves on $X$.
Assume that for any $s>0$, there is $N(s) > 0$ such that for any $n > N(s)$, $(L^s)^0 E_n \in V_{\coh(X)}^{\Phi}(3)$.
Then
$\mu^+(E_n) \to 0$ as $n \to +\infty$.
\end{prop}
\begin{proof}
Assume $\mu^+(E_n) \not \to 0$ as $n \to +\infty$ for a contradiction.
Then there exists $\varepsilon > 0$ such that
for any $N'>0$ there is $n>N'$ satisfying $\mu^+(E_n) > \varepsilon$.

Let $T_n$ be the slope semistable H-N factor of $E_n$ with the highest slope, i.e. $\mu(T_n) = \mu^+(E_n)$.
There is $s \in \N$ such that $\mu ((L^s)^0) > -\varepsilon$. Then for any $N'>0$ there is $n>N'$ such that
$(L^s)^0 T_n \in \T_0$.
Therefore, for some $n > N(s)$ we have
\begin{align*}
\Hom((L^s)^0 T_n, (L^s)^0 E_n) & \cong \Hom(\Phi((L^s)^0 T_n), \Phi((L^s)^0 E_n)) \\
                               & \cong 0,
\end{align*}
as $((L^s)^0 T_n)^3 =0$ (from Proposition \ref{prop4.3}) and $(L^s)^0 E_n \in V_{\coh(X)}^{\Phi}(3)$.
This is the required contradiction to complete the proof.
\end{proof}
Let $s$ be a positive integer. Consider the Fourier-Mukai functor defined by
$$
\Pi = \Phi \circ (L^s)^0 \circ \Phi [3].
$$
Then $\Pi^i_{\coh(X)}(\OO_x) = 0$ for $i \ne 0$ and $\Pi^0_{\coh(X)}(\OO_x) = L^s  \Po_y$ for some $y \in X$.
Define the Fourier-Mukai functor
$$
\widehat{\Pi} = \Phi \circ (L^{-s})^3 \circ \Phi.
$$
One can show that $\widehat{\Pi}[3]$ is right and left adjoint to $\Pi$ (and vice versa).
We have  $\widehat{\Pi}^i_{\coh(X)}(\OO_x) = 0$ for $i \ne 0$, and
$\widehat{\Pi}^0_{\coh(X)}(\OO_x) = L^{-s}  \Po_z$ for some $z \in X$.
Therefore $\Pi$ is a Fourier-Mukai functor with kernel a locally free sheaf $\mathcal{U}$ on $X \times X$.

We have the spectral sequence
\begin{equation}
\label{auxss}
\Phi^p \left( (L^s)^0 \, \Phi^q(E) \right) \Longrightarrow \Pi^{p+q-3}(E)
\tag{\ding{71}}
\end{equation}
for $E$.
\begin{prop}
\label{prop4.9.2}
Let $E$ be a coherent sheaf such that $L^{-n}E \in V_{\coh(X)}^{\Phi}(k)$ for sufficiently large $n$, where $k\in \{0, \ldots, 3\}$.
Then $\mu^+((L^{-n}E)^k) \to 0$ as $n \to +\infty$.
\end{prop}
\begin{proof}
Since $L^{-n}E \in V_{\coh(X)}^{\Phi}(k)$ for sufficiently large $n$,  $E \in \coh^k(X)$.
If $k=0$ then $E \in  \coh^0(X)$ and so we have $\mu^+((L^{-n}E)^0) =0$.
Otherwise, by Propositions \ref{prop4.8}, \ref{prop4.7} and \ref{prop4.4}, for $E \in \coh^k(X)$ we have
$(L^{-n}E)^k \in \T_0$.
Let $s$ be a positive integer. Consider the convergence of the Spectral Sequence \eqref{auxss}.
For large enough $n$, we also have $L^{-n}E \in V_{\coh(X)}^{\Pi}(k)$.
Therefore $(L^s)^0  (L^{-n}E)^k \in V_{\coh(X)}^{\Phi}(3)$. By Proposition \ref{prop4.9.1} we have
$\mu^+((L^{-n}E)^k) \to 0$ as $n \to +\infty$.
\end{proof}
\begin{prop}
\label{prop4.9.3}
Let $E$ be  a reflexive sheaf. Then for sufficiently large $n > 0$,
\begin{enumerate}[label=(\roman*)]
\item $L^{-n}E \in V^{\Phi}_{\coh(X)}(2,3)$, and
\item $(L^{-n}E)^2 \cong (T_0)^0$ for some $T_0 \in \coh^0(X)$.
\end{enumerate}
\end{prop}
\begin{proof}
\begin{enumerate}[label=(\roman*)]
\item
Consider  a minimal locally free resolution of $E$, 
$$
0 \to F_2 \to F_1 \to E \to 0.
$$
By applying the FMT $\Phi L^{-n} $ for  sufficiently large  $n > 0$, we obtain
$L^{-n}E \in V^{\Phi}_{\coh(X)}(2,3)$.
\item
Since $E$ is a reflexive sheaf, there is a locally free sheaf $P$ and a torsion free sheaf $Q$ such that
$$
0 \to E \to P \to Q \to 0
$$
is a SES in $\coh(X)$.
By applying the FMT $\Phi L^{-n}$ for sufficiently large $n$ we have $(L^{-n}E)^2 \cong (L^{-n}Q)^1$.

The torsion free sheaf $Q$ fits into the SES $0 \to Q \to Q^{**} \to T \to 0$ for
some $T \in  \coh^{\le 1}(X)$. Apply the FMT $\Phi L^{-n}$ for sufficiently large $n$ and consider the LES of $\coh(X)$-cohomologies. Since
$L^{-n}Q^{**}\in V_{\coh(X)}^{\Phi}(2,3)$, we have $(L^{-n}Q)^1 \cong (L^{-n}T)^0$.
The torsion sheaf $T \in \coh^{\le 1}(X)$ fits into a SES
$0 \to T_0 \to T \to T_1 \to 0$ in $\coh(X)$ for $T_i \in \coh^i(X)$, $i=0,1$.
Therefore   $(L^{-n}T)^0 \cong  (T_0)^0$, and so $(L^{-n}E)^2 \cong (T_0)^0$ as required.
\end{enumerate}
\end{proof}
\begin{prop}
\label{prop4.9.4}
Let $E \in \coh^1(X)$ with $E \in V^{\Phi}_{\coh(X)}(1)$.
If $0 \ne T \in \HN[0,+\infty]$ is a subsheaf of $E^1$, then
$\ell \ch_2(T) \le 0$.
\end{prop}
\begin{proof}
For $n>0$ and generic $z \in X$, we have the $\coh(X)$-SES
$$
0 \to L^{-n}_z E  \to E  \to T_0 \to 0
$$
for $T_0 :=  E|_{ nD_z} \in \coh^0(X)$.
By applying the FMT $\Phi$, we get the following commutative diagram
$$
\xymatrixcolsep{2pc}
\xymatrixrowsep{2pc}
\xymatrix{
0 \ar[r]  & T_0^0 \ar[r]   &  ( L^{-n}_z E)^1 \ar[r]  &   E^1 \ar[r] & 0 \\
0 \ar[r] & T_0^0 \ar[r] \ar@{=}[u] & A \ar[r] \ar@{^{(}->}[u] &   T \ar[r] \ar@{^{(}->}[u] & 0,
}
$$
for some $A \in \HN[0,+\infty]$. 
So we have $\ch_k(A) = \ch_k(T)$ for $k=1,2,3$.

Let $G$ be a slope semistable H-N factor of $A$. From the usual B-G inequality,
 $\ell (\ch_1(G)^2 - 2 \ch_0(G) \ch_2(G)) \ge 0$. So we have
\begin{align*}
2 \ell \ch_2(G) & \le \frac{\ell \ch_1(G)^2}{\ch_0(G)} \\
                & = \ell^2 \ch_1(G) \, \mu(G) \\
                & \le \ell^2 \ch_1(A) \, \mu(G) \\
                & \le  \ell^2 \ch_1(T) \,  \mu^+ \left(( L^{-n}_z E)^1 \right).
\end{align*}
By Proposition \ref{prop4.9.2}, $\mu^+ (( L^{-n}_z E)^1) \to 0$ as $n \to +\infty$.
So choose large enough $n>0$ such that
$ \ell^2 \ch_1(T)  \mu^+ (( L^{-n}_z E)^1) <\ell^3$.
Since $2 \ell \ch_2(G) \in \ell^3 \Z$ we have $\ell \ch_2(G) \le 0$.
So $\ell \ch_2(T)= \ell \ch_2(A) \le 0$.
\end{proof}
\begin{prop}
\label{prop4.9.5}
We have the following:
\begin{enumerate}[label=(\roman*)]
\item
Let $E \in \F_0$  be a reflexive sheaf.
If $0 \ne T \in \T_0$ is a subsheaf of $E^1$ then $\ell \ch_2(T) \le  0$.

\item
Let  $E \in \T_0$  be a torsion free.
If  $0 \ne F \in \F_0$ is a quotient of $E^2$ then $\ell \ch_2(F) \le  0$.
\end{enumerate}
\end{prop}
\begin{proof}
\begin{enumerate}[label=(\roman*)]
\item

Recall that, for any positive integer $m$, non-reduced divisors $mD_x$ of $L_x^m$ are 
topologically supported on $D_x$.

Since $E$ is a reflexive sheaf, one can choose $x, y \in X$ such that
\begin{itemize}
\item $\dim (D_x \cap D_y) = 1$,
\item $E|_{D_x}$  is locally free on $D_x$, and
\item $E|_{D_y}$ is locally free on $D_y$.
\end{itemize}

By Proposition \ref{prop4.9.3}, for sufficiently large $m>0$, $L_x^{-m} E \in V^{\Phi}_{\coh(X)}(2,3)$.
By applying the FMT $\Phi$ to the $\coh(X)$-SES
$$
0 \to L_x^{-m} E \to E \to E|_{mD_x} \to 0
$$
$E|_{mD_x} \in V^{\Phi}_{\coh(X)}(1,2)$ and $E^1 \hookrightarrow \left(E|_{mD_x} \right)^1$.
Since $E|_{D_x}$  is locally free on $D_x$, for large enough $n>0$,   $L_y^{-n} E|_{mD_x} \in
V^{\Phi}_{\coh(X)}(2)$ .
By applying the FMT $\Phi$ to the $\coh(X)$-SES
$$
0 \to L_y^{-n} E|_{mD_x} \to E|_{mD_x} \to E|_{mD_x \cap nD_y} \to 0,
$$
$E|_{mD_x \cap nD_y} \in V^{\Phi}_{\coh(X)}(1)$ and $\left(E|_{mD_x} \right)^1 \hookrightarrow \left( E|_{mD_x \cap nD_y} \right)^1$.
Therefore we have
$$
T \hookrightarrow E^1 \hookrightarrow \left(  E|_{mD_x \cap nD_y} \right)^1.
$$
The result follows from Proposition \ref{prop4.9.4}.

\item
Since $F \ne 0$ is a quotient of $E^2$, we have
$F^{*} \hookrightarrow (E^2)^*$.
Here $F^* \in \HN[0, +\infty)$ fits into $\coh(X)$-SES
$0 \to T \to F^* \to F_0 \to 0$ for some $T \in \T_0$ and $F_0 \in \HN[0]$.
By the usual B-G inequality $\ell \ch_2(F_0) \le 0$.

By Proposition \ref{prop4.3}, $E^3 = 0 = (E^*)^0$.
Therefore from the convergence of the ``Duality'' Spectral Sequence \ref{dualss}
for $E$, we have the $\coh(X)$-SES
$$
0 \to (-1)^* (E^*)^1 \to (E^2)^*  \to P \to 0,
$$
for some subsheaf $P$ of $(\calExt^1(E, \OO))^0$.
By Proposition \ref{prop4.4}, $(\calExt^1(E, \OO))^0 \in \F_0$ and so $P \in \F_0$.
Therefore $\Hom (T, P) = 0$ and so $P \hookrightarrow  (-1)^* (E^*)^1$.
Here $E^* \in \F_0$ and so by part (I), $\ell \ch_2(T) \le 0$.
Therefore  $\ell \ch_2(F) \le \ell \ch_2(F^{**}) = \ell \ch_2(F^{*}) = \ell \ch_2(F_0) + \ell \ch_2(T) \le 0$.

\end{enumerate}
\end{proof}
\begin{prop}
\label{prop4.12}
For $E \in \coh(X)$
\begin{enumerate}[label=(\roman*)]
\item if $E \in \F_0$ then  $E^1 \in \F_0$, and
\item if $E \in \HN[0,+\infty)$ with $E^3 = 0$ then  $E^2 \in \HN[0,+\infty]$.
\end{enumerate}
\end{prop}
\begin{proof}
\begin{enumerate}[label=(\roman*)]
\item
Assume the opposite for a contradiction.
Let $T = T(E^1)$ and $F= F(E^1)$. Then $0 \to T \to E^1 \to F \to 0$ is a  SES in $\coh(X)$.
By Proposition \ref{prop4.6} $E^1$ is reflexive and so non-trivial $T$ is reflexive.
So $\ell^2 \ch_1(T) > 0$.
By applying the FMT $\Phi$ to this SES we obtain that
$T \in V_{\coh(X)}^{\Phi}(2)$ and $F \in V_{\coh(X)}^{\Phi}(1,2)$.
Moreover, we have the $\coh(X)$-SES
$$
0 \to F^1 \to T^2 \to E_1 \to 0
$$
for some subsheaf $E_1$ of $E^{12}$.
From the Mukai Spectral Sequence \ref{mukaiss}
for $E$ we have the $\coh(X)$-SES
$$
0 \to E^{20} \to E^{12} \to E_2 \to 0,
$$
for some subsheaf $E_2$ of $(-1)^* E$. Therefore $E_2 \in \F_0$ and by Proposition \ref{prop4.4}
$E^{20} \in  \F_0$. So we have $E^{12}   \in \F_0$. Hence  $E_1 \in \F_0$.

Let $T_1 : = T(F^1)$ and $F_1 : = F(T^2)$.
They fit into the following commutative diagram for some $F_2 \in \F_0$.
$$
\xymatrixcolsep{2pc}
\xymatrixrowsep{1.5pc}
\xymatrix{
 & 0 & 0 &   &  \\
 0 \ar[r] & F_2 \ar[r]\ar[u] & F_1 \ar[r]\ar[u] & E_1 \ar[r] & 0 \\
 0 \ar[r] & F^1 \ar[r] \ar[u] &  T^2 \ar[r] \ar[u] &  E_1 \ar[r] \ar@{=}[u] & 0 \\
    & T_1 \ar[u] \ar@{=}[r]  & T_1 \ar[u] &   &  \\
     & 0 \ar[u] & 0 \ar[u] &   &   }
$$
By Proposition \ref{prop4.9.5},   $\ell \ch_2(F_1) \le 0$.

By applying the FMT $\Phi$ to the $\coh(X)$-SES
$0 \to T_1 \to F^{1} \to F_2 \to 0$ we obtain the $\coh(X)$-SES
$$
0 \to F_2^1 \to T_1^2 \to F_3 \to 0
$$
for some subsheaf $F_3$ of $F^{12}$. Also
$T_1 \in V^{\Phi}_{\coh(X)}(2)$.
By considering the Mukai Spectral Sequence \ref{mukaiss}
for $F$, one can show $F_3 \in \F_0$.
By Proposition \ref{prop4.6} $F_2^1$ is reflexive. So
$T_1^2$ is torsion free and it fits into $\coh(X)$-SES
$$
0 \to T_1^2 \to (T_1^2)^{**} \to Q \to 0,
$$
for some $Q \in \coh^{\le 1}(X)$. The torsion sheaf $Q$ fits into $\coh(X)$-SES
$$
0 \to Q_0 \to Q \to Q_1 \to 0
$$
for $Q_0 \in \coh^0(X)$ and $Q_1 \in \coh^1(X)$.
By Proposition \ref{prop4.9.3}, for large enough $m>0$, $(L^{-m}T_1^2)^1 \cong (L^{-m}Q)^0 \cong Q_0^0$.
Also $(L^{-m}Q_1)^1 \cong (L^{-m}Q)^1$  and $(L^{-m}(T_1^2)^{**})^2 \cong R_0^0$ for some $R_0 \in \coh^0(X)$.
So we have the $\coh(X)$-SES
$$
0 \to (L^{-m}Q_1)^1 \to (L^{-m}T_1^2)^2 \to R_0^0 \to 0.
$$
By Proposition \ref{prop4.8}, $(L^{-m}T_1^2)^2 \in \HN[0, +\infty)$, and $\ell \ch_2((L^{-m}T_1^2)^2) = 0$.

The torsion free sheaf $F_3$ also fits into $\coh(X)$-SES
$0 \to F_3 \to F_3^{**} \to S \to 0$ for some
$S \in \coh^{\le 1}(X)$.

Choose $x, y \in X$ such that
\begin{itemize}
\item $\dim(D_x \cap D_y) =1$,
\item $D_{x}  \cap \Supp(Q_0) = \emptyset$,
\item $\dim(\Supp(Q_1) \cap D_x) \le 0$,
\item $D_{x}  \cap D_{y}  \cap \Supp(Q) = \emptyset$,
\item $D_{x}  \cap D_{y}  \cap \Supp(S) = \emptyset$,
\item since $F_2^1$ is reflexive, $F_2^1|_{D_x}$  is locally free on $D_x$, and $F_2^1|_{D_y}$  is locally free on $D_y$,
\item since $F_3^{**}$ is reflexive, $F_3^{**}|_{D_x}$  is locally free on $D_x$, and $F_3^{**}|_{D_y}$  is locally free on $D_y$.
\end{itemize}

From the Mukai Spectral Sequence for $F_2$, $F_2^1 \in V^{\Phi}_{\coh(X)}(2,3)$.
Since it is a reflexive sheaf,  for large enough $m>0$, $L^{-m}_xF_2^1 \in V^{\Phi}_{\coh(X)}(2,3)$, and
since $F_2^1|_{D_x}$  is locally free on $D_x$, $L^{-n}_yF_2^1|_{mD_x} \in V^{\Phi}_{\coh(X)}(2)$.
So $F_2^1|_{mD_x \cap nD_y} \in V^{\Phi}_{\coh(X)}(1)$.
Since $D_{x}  \cap D_{y} \cap \Supp(S) = \emptyset$, similarly one can show
$F_3 |_{mD_x \cap nD_y} \cong F_3^{**}|_{mD_x \cap nD_y} \in V^{\Phi}_{\coh(X)}(1)$.
Therefore we have $ T_1^2|_{mD_x \cap nD_y} \cong  (T_1^2)^{**}|_{mD_x \cap nD_y} \in V^{\Phi}_{\coh(X)}(1)$.

By applying the FMT $\Phi$ to the $\coh(X)$-SES
$ 0 \to L_x^{-m}T_1^2 \to T_1^2 \to T_1^2|_{mD_x} \to 0$, for large enough $m>0$
we have the $\coh(X)$-LES
$$
0 \to Q_0^0 \to (-1)^* T_1 \to (T_1^2|_{mD_x})^1 \to (L^{-m}T_1^2)^2 \to 0.
$$
So $\left(T_1^2|_{mD_x}\right)^1 \in \HN[0, +\infty]$ and
$\ch_2((T_1^2|_{mD_x})^1) = \ch_2(T_1)$.
Moreover we have the $\coh(X)$-SES
$$
0 \to T_1^2|_{mD_x} \to (T_1^2)^{**}|_{mD_x} \to Q_1|_{mD_x} \to 0.
$$
Here $Q_1|_{mD_x} \in \coh^0(X)$.
So for   large enough $n>0$, $(L^{-n}T_1^2|_{mD_x})^1 \cong (Q_1|_{mD_x})^0$.

By applying the FMT $\Phi$ to the $\coh(X)$-SES
$0 \to L_y^{-n} T_1^2|_{mD_x} \to T_1^2|_{mD_x} \to T_1^2|_{mD_x \cap nD_y} \to 0$
we have the $\coh(X)$-LES
$$
0 \to (Q_1|_{mD_x})^0 \stackrel{\alpha}{\to} \left(T_1^2|_{mD_x}\right)^1 \to (T_1^2|_{mD_x \cap nD_y})^1 \to \cdots.
$$
Let $T_2 : = \coker(\alpha)$. Then $T_2 \in \HN[0, +\infty]$ and $\ch_2(T_2) = \ch_2(T_1)$.
By Proposition \ref{prop4.9.4}, we have $\ell \ch_2(T_2) \le 0$.
So  $\ell \ch_2(T^2) = \ell \ch_2(T_1) + \ell \ch_2(F_1) \le 0$.
Therefore we have
$\ell^2 \ch_1(T) \le 0$. This is the required contradiction to complete the proof.

\item
Since $E^* \in \HN(-\infty,0]$, from (i)  $(E^*)^1 \in \HN(-\infty,0]$.
By the convergence of the ``Duality'' Spectral Sequence \ref{dualss}
for $E$ we have $(E^2)^* \in  \HN(-\infty,0]$.
So $E^2 \in  \HN[0,+\infty]$ as required.
\end{enumerate}
\end{proof}
\begin{cor}
\label{prop4.13}
Let $E \in \T_0$. Then $E^2 \in \T_0$.
\end{cor}
\begin{proof}
Let $T = T(E^2)$ and $F= F(E^2)$. Then $0 \to T \to E^2 \to F \to 0$ is a SES in $\coh(X)$.
Now we need to show $F = 0$.
Apply the FMT  $\Phi$ and consider the LES of $\coh(X)$-cohomologies. So we have $F \in V^{\Phi}_{\coh(X)}(1)$ and
$$
0 \to T^1 \to E^{21} \to F^1 \to T^2 \to 0
$$
is a LES in $\coh(X)$. From the convergence of the Mukai Spectral
Sequence \ref{mukaiss} for $E$ we have the $\coh(X)$-SES
$$
0 \to Q \to E^{21} \to E^{13} \to 0,
$$
where $Q$ is a quotient of $(-1)^* E$. Then $Q \in\T_0$ and, by
Proposition \ref{prop4.4}, $E^{13}\in\T_0$ and so $ E^{21}  \in \T_0$.
On the other hand, by Proposition \ref{prop4.12}, $F^1 \in \F_0$.
So the map $E^{21} \to F^1$ is zero and  $F^1 \cong T^2$. Hence
 $F \cong (-1)^* F^{12} \cong (-1)^* T^{22} = 0$ (from the Mukai
 Spectral Sequence \ref{mukaiss} for $T$) as required.
\end{proof}
\begin{prop}
\label{prop4.15}
Let $E \in \HN(0,1]$. Then $E^0 \in \HN(-\infty,-\frac{1}{2}]$.
\end{prop}
\begin{proof}
Due to Mukai, $\Phi L \Phi \cong (-1)^* L^{-1} \Phi L^{-1}$. Therefore
we have the following convergence of spectral sequence:
$$
E_2^{p,q} = \Phi^p_{\coh(X)} L \Phi^q_{\coh(X)}(E) \Longrightarrow (-1)^* L^{-1} \Phi^{p+q}_{\coh(X)} (L^{-1}E).
$$
Here $L^{-1}E \in \HN(-1,0]$, and so by Proposition
  \ref{prop4.3}, $(L^{-1}E)^0 = 0$. So from the convergence of the above spectral
  sequence for $E$ we have $(LE^0)^0 = 0$.
  Also $(LE^0)^1 \hookrightarrow L^{-1}(L^{-1}E)^1$. By Proposition \ref{prop4.12} $(L^{-1}E)^1 \in \F_0$ and so
  $(LE^0)^1 \in \HN(-\infty,-1] \subset \F_0$.

Let $F \subset E^0$ be the H-N semistable factor of $E^0$ with the
highest slope and let $\mu: = \mu(F)$. Then
$(LF)^0\hookrightarrow (LE^0)^0$ and so $(LF)^0=0$. Let $\ch(F) = (a_0, \mu a_0, a_2, a_3)$. Now suppose
$\mu > -\frac{1}{2}$ for a contradiction. Then $LF\in\T_0$ and $F$ fits into the $\coh(X)$-SES
\begin{equation}
\label{SES-mu-half}
0 \to F \to E^0 \to G \to 0,
\tag{\ding{73}}
\end{equation}
for some $G \in \HN(-\infty, 0]$. By Proposition \ref{prop4.2}, $E^0$
is reflexive. Since $G$ is torsion-free, it follows that $F$ is also reflexive.
Apply the FMT $\Phi$ and consider the LES of $\coh(X)$-cohomologies. Then  we have $F \in V^{\Phi}_{\coh(X)}(2,3)$ and
$$
0 \to G^1 \to F^2 \to E^{02} \to \cdots
$$
is an exact sequence in $\coh(X)$.
From the convergence of the Mukai Spectral Sequence \ref{mukaiss} for
$E$,  $E^{02}\cong E^{10}$ and $E^{10}\in H(-\infty, 0]$ by Proposition
  \ref{prop4.4}. Also by Proposition \ref{prop4.12}, $G^1 \in
  \HN(-\infty, 0]$. So
$F^2 \in \HN(-\infty, 0]$ and we have $\ell^2 \ch_1(F^2) \le 0$. Moreover, by Proposition \ref{prop4.4}, $F^3 \in \HN(0,
+\infty]$ and so
$\ell^2 \ch_1(F^3) \ge 0$. Therefore $\ell^2 \ch_1(\Phi(F)) \le 0$ and so $\ch(\Phi(F))=(a_3, -a_2, \mu a_0, -a_0)$ implies
$$
a_2 \ell^3 = 2 \ell \ch_2(F) \ge 0.
$$
Apply the FMT $\Phi L$ to the SES \eqref{SES-mu-half} and consider the LES of $\coh(X)$-cohomologies.
Then we have the $\coh(X)$-LES
$$
0 \to (LG)^0 \to (LF)^1 \to (LE^0)^1 \to \cdots.
$$
Here $ (LE^0)^1 \in \F_0$ and so $(LF)^1 \in \F_0$.
By Corollary \ref{prop4.13}  $(LF)^2 \in \HN(0, +\infty]$.
So $ \ell^2 \ch_1(LF^1) \le 0$ and $ \ell^2 \ch_1(LF^2) \ge 0$ which imply
$\ell^2 \ch_1(\Phi(LF)) \ge 0$. Hence
$$
(a_0 + 2 \mu a_0 + a_2) \ell^3 = 2\ell \ch_2(LF) \le 0.
$$
Here by the assumption $2\mu + 1 > 0$ and we already obtained that $a_2  \ge 0$. Hence $(2 \mu+1)a_0 + a_2 > 0$ and which is not possible. 
This is the required contradiction to complete the proof.
\end{proof}
\begin{prop}
\label{prop4.16}
Let $E \in \HN[-1,0]$. Then $E^3 \in \HN[\frac{1}{2}, +\infty]$.
\end{prop}
\begin{proof}
From the ``Duality'' Spectral Sequence \ref{dualss} for $E$ we have $( E^* )^0\cong (-1)^*(E^3)^*$.
Here $E^{*} \in \HN[0,1]$ and so by Propositions \ref{prop4.3} and \ref{prop4.15}, $(E^*)^0 \in \HN(-\infty, -\frac{1}{2}]$.
Hence $(E^3)^* \in \HN(-\infty, -\frac{1}{2}]$ and so
$E^3 \in \HN[\frac{1}{2}, +\infty]$ as required.
\end{proof}
\begin{thm}
\label{prop4.17}
We have the following:
\begin{enumerate}[label=(\roman*)]
\item $L \Phi\left( \B \right)  \subset \langle \B,\B[-1], \B[-2] \rangle$, and
\item $\Phi L^{-1} [1] \left( \B \right)  \subset \langle \B,\B[-1], \B[-2] \rangle$.
\end{enumerate}
\end{thm}
\begin{proof}
\begin{enumerate}[label=(\roman*)]
\item We can visualize $\B$ as follows:
$$
\begin{tikzpicture}[scale=1.2]
\draw[style=dashed] (0.5,0) grid (3.5,1);
\fill[lightgray] (1,0) -- (2,0) -- (3,1) -- (2,1) --(1,0);
\draw[style=thick] (1,0) -- (2,0)-- (3,1)-- (2,1)-- (1,0);
\draw[style=thick] (2,0) -- (2,1);
\draw (1.7,0.3) node {$\scriptstyle B$};
\draw (2.3,0.7) node {$\scriptstyle A$};
\draw (1.5,- 0.3) node {$\scriptstyle{-1}$};
\draw (2.5,-0.3) node {$\scriptstyle{0}$};
\draw (-0.7,0.5) node {$\mathcal{B} = \langle \mathcal{F}[1] , \mathcal{T} \rangle : $};
\draw (6,0.5) node {$\scriptstyle A \in \mathcal{T} = \HN(\frac{1}{2}, +\infty], \ B \in \mathcal{F} = \HN(-\infty, \frac{1}{2}]$ };
\end{tikzpicture}
$$
If $E \in \F = \HN(-\infty, \frac{1}{2}]$ then by Propositions \ref{prop4.3} and \ref{prop4.15}, $L E^0 \in  \F$. Also by Proposition \ref{prop4.4},
$LE^3 \in \HN(1, +\infty] \subset \F$. Therefore $ L \Phi(E)$ has $\B$-cohomologies in 1,2,3 positions. That is
$$
L \Phi \left( \F \right)[1] \subset \langle \B,\B[-1], \B[-2] \rangle.
$$
$$
\begin{tikzpicture}[yscale=1.2]
\draw[style=dashed] (0.5,-2) grid (6.5,-1);
\fill[lightgray] (1,-2) -- (4,-2) -- (5,-1) -- (2,-1) --(1,-2);
\draw[style=thick] (1,-2) -- (4,-2) -- (5,-1) -- (2,-1) --(1,-2);
\draw[style=dashed] (2,-2) -- (2,-1) ;
\draw[style=dashed] (3,-2) -- (3,-1) ;
\draw[style=dashed] (4,-2) -- (4,-1) ;
\draw (1.65,-1.75) node {$\scriptscriptstyle LB^0$};
\draw (2.5,-1.5) node {$\scriptscriptstyle LB^1$};
\draw (3.5,-1.5) node {$\scriptscriptstyle LB^2$};
\draw (4.35,-1.25) node {$\scriptscriptstyle LB^3$};
\draw[style=thick] (2,-2) -- (2.3,-1.7);
\draw[style=thick] (2.7,-1.3) -- (3,-1);
\draw[style=thick] (3,-2) -- (3.3,-1.7);
\draw[style=thick] (3.7,-1.3) -- (4,-1);
\draw (1.5,- 2.3) node {$\scriptscriptstyle{-1}$};
\draw (2.5,-2.3) node {$\scriptscriptstyle{0}$};
\draw (3.5,- 2.3) node {$\scriptscriptstyle{1}$};
\draw (4.5,-2.3) node {$\scriptscriptstyle{2}$};
\draw (5.5,-2.3) node {$\scriptscriptstyle{3}$};
\draw[style=dashed] (-3 ,-2) grid (-1,-1);
\fill[lightgray] (-3,-2) -- (-2,-2) -- (-2,-1) -- (-3,-2);
\draw[style=thick] (-3,-2) -- (-2, -2)-- (-1, -1) -- (-2, -1) -- (-3, -2);
\draw[style=thick] (-2,-2) -- (-2, -1);
\draw (-2.3,-1.7) node {$\scriptstyle B$};
\draw (-2.5,- 2.3) node {$\scriptscriptstyle{-1}$};
\draw (-1.5,-2.3) node {$\scriptscriptstyle{0}$};
\draw (-3.2,-2.1) to [out=120,in=240] (-3.2,-0.9);
\draw (-0.8,-2.1) to [out=60,in=300] (-0.8,-0.9);
\draw (-3.7,-1.5) node {${L \Phi}$};
\draw (0,-1.5) node {$=$};
\end{tikzpicture}
$$

On the other hand if $E \in \T = \HN(\frac{1}{2}, +\infty]$ then by Proposition \ref{prop4.3} $L E^3 =0$ and by Corollary \ref{prop4.13}  $LE^2 \in \HN(1, +\infty] \subset \T$.
So $L\Phi(E)$ has $\B$-cohomologies in 0,1,2 positions. That is
$$
L \Phi\left( \T \right) \subset \langle \B,\B[-1], \B[-2] \rangle.
$$
$$
\begin{tikzpicture}[yscale=1.2]
\draw[style=dashed] (0.5,0) grid (6.5,1);
\fill[lightgray] (2,0) -- (4,0) -- (5,1) -- (2,1) --(2,0);
\draw[style=thick] (2,0) -- (4,0)  -- (5,1) -- (2,1) -- (2,0);
\draw[style=dashed] (3,0) -- (3,1) ;
\draw[style=dashed] (4,0) -- (4,1) ;
\draw (2.5,0.5) node {$\scriptscriptstyle LA^0$};
\draw (3.5,0.5) node {$\scriptscriptstyle LA^1$};
\draw (4.35,0.75) node {$\scriptscriptstyle LA^2$};
\draw[style=thick] (2,0) -- (2.3,0.3);
\draw[style=thick] (2.7,.7) -- (3,1);
\draw[style=thick] (3,0) -- (3.3,0.3);
\draw[style=thick] (3.7,.7) -- (4,1);
\draw (1.5,- 0.3) node {$\scriptscriptstyle{-1}$};
\draw (2.5,-0.3) node {$\scriptscriptstyle{0}$};
\draw (3.5,- 0.3) node {$\scriptscriptstyle{1}$};
\draw (4.5,-0.3) node {$\scriptscriptstyle{2}$};
\draw (5.5,-0.3) node {$\scriptscriptstyle{3}$};
\draw[style=dashed] (-3 ,0) grid (-1,1);
\fill[lightgray] (-2,0) -- (-1,1) -- (-2,1) -- (-2,0);
\draw[style=thick] (-3,0) -- (-2, 0)-- (-1, 1) -- (-2, 1) -- (-3, 0);
\draw[style=thick] (-2,0) -- (-2, 1);
\draw (-1.7,0.7) node {$\scriptstyle A$};
\draw (-2.5,- 0.3) node {$\scriptscriptstyle{-1}$};
\draw (-1.5,-0.3) node {$\scriptscriptstyle{0}$};
\draw (-3.2,-0.1) to [out=120,in=240] (-3.2,1.1);
\draw (-0.8,-0.1) to [out=60,in=300] (-0.8,1.1);
\draw (-3.7,0.5) node {${L \Phi}$};
\draw (0,0.5) node {$=$};
\end{tikzpicture}
$$

Hence $L \Phi\left( \B \right) \subset \langle \B,\B[-1], \B[-2] \rangle$,
as $\B = \langle \F[1] ,\T \rangle$.

\item We can use Propositions \ref{prop4.3}, \ref{prop4.4} and
  \ref{prop4.16}, and Proposition \ref{prop4.12} in a  similar way to the proof of (i).
\end{enumerate}
\end{proof}
\section{(Semi)stable sheaves with the Chern character $(r, 0, 0, \chi)$}
\label{section5}
\label{sect1}

In this section we shall consider sheaves $E$ with
$\ch_k(E)=0$ for $k=1,2$ which arise as
the $\coh(X)$-cohomology of some of
the tilt-stable objects.
For example, when $F \in \B$ is a tilt stable object with $\nu(F) =0$ and $F_i :=
H^{i}_{\coh(X)}(F)$. By Proposition \ref{prop3.1},
if $\mu(F_{-1})=0$ then $\ch_k(F_{-1})=0$, and if $\mu(F_{0})=1$ then $\ch_k(L^{-1}F_{0})=0$ for $k=1$ and $2$.

We would like to show that such sheaves can only take a very special form:
\begin{thm}
\label{prop5.1}
Let $E$ be a slope semistable sheaf with $\ch_k(E)=0$ for $k=1,2$. Then
$E^{**}$ is a homogeneous bundle. In other words, $E^{**}$ is filtered with quotients from $\Pic^0(X)$.
\end{thm}
\begin{proof}
Assume the opposite for a contradiction.
Then there exists a stable reflexive sheaf $E$ with $\ch_k(E) = 0$ for $k=1,2$,  and $H^k(X, E \otimes \Po_x) = 0$ for $k=0,3$ and any $x \in X$.
By a result of Simpson (\cite[Theorem 2]{Sim}) we have $\ch_3(E) = 0$.
Therefore, $\ch(E) = (r, 0 , 0, 0)$ for some positive integer $r$.

Since $H^k(X, E \otimes \Po_x) = 0$ for $k=0,3$ and any $x \in X$, we have $E^0 = E^3 =0$.
By Proposition \ref{prop4.12},  $E^1 \in \HN(-\infty, 0]$ and  $E^2 \in \HN[0,+\infty]$.
So we have $\ell^2 \ch_1(E^1) \le 0$ and $ \ell^2 \ch_1(E^2) \ge
0$. Therefore, $\ell^2 \ch_1(\Phi(E)) \ge 0$ which implies $\ell \ch_2(E) \le
0$. Since $\ch_2(E) =0$, we obtain $\ch_1(E^1) = \ch_1(E^2) =0$.
Then we have
\begin{align*}
\ch(E^1) = (a, 0 , -b, c),  \ \ \ch(E^2) = (a, 0, -b, -r+c),
\end{align*}
for some $a>0$ and $b \ge 0$. Moreover we have $E^1 \in \HN[0]$.

If $E^{13} \ne 0$ then $E^1$ fits into a $\coh(X)$-SES of the form
$
0 \to K_1 \to E^1 \to \Po_{z_1} \I_{C_1} \to 0.
$
Then $K_1 \in \HN[0]$ and we have the following exact sequence
$$
\cdots \to K_1^3 \to E^{13} \to \OO_{-z_1} \to 0
$$
in $\coh(X)$.
If $K_1^3 \ne 0$ then $K_1$  fits into a $\coh(X)$-SES
$
0 \to K_2 \to K_1 \to \Po_{z_2} \I_{C_2} \to 0.
$
Then $K_2 \in \HN[0]$ and we have the following exact sequence
$$
\cdots \to K_2^3 \to K_1^3 \to \OO_{-z_2} \to 0
$$
in $\coh(X)$.
We can continue this process for only a finite number of steps since $\rk (E^1) < +\infty$ and hence $E^{13}$ is filtered by skyscraper sheaves.
Moreover from the convergence of the Mukai Spectral Sequence \ref{mukaiss} for $E$, we have the $\coh(X)$-SES
$$
0 \to E^{20} \to E^{12} \to  Q \to 0
$$
where $Q$ is a subsheaf of $(-1)^*E$ and so $Q \in \HN(-\infty,
0]$. By Proposition \ref{prop4.4}, $E^{20}\in\HN(-\infty,0]$. This implies $E^{12} \in \HN(-\infty, 0]$.
Then $\ell^2 \ch_1(\Phi(E^1)) \le 0$ and so $-b \ell^3 =  2 \ell \ch_2(E^1) \ge 0$. Hence $b = 0$.
By Proposition \ref{prop4.6}, $E^1$ is a reflexive sheaf and since $E^1 \in \HN[0]$ it is slope semistable.
So by \cite[Theroem 2]{Sim} we have $c= \ch_3(E^1) =0$.
Therefore $\ch(\Phi(E^1)) = (0,0,0,-a)$. Since $E^{13} \in \coh^0(X)$, we have $\ch_k(E^{12}) = 0$ for $k=0,1,2$.
So $E^{12} \in \HN(0, +\infty]$. Therefore $E^{12} = 0$ and we have the $\coh(X)$-SES
$$
0 \to (-1)^*E \to E^{21} \to E^{13} \to 0.
$$
Since $E^{13} \in \coh^0(X)$ and $E$ is locally free, $\Ext^1(E^{13}, (-1)^* E) = 0$.
Therefore $E^{21}\cong (-1)^*E \oplus E^{13}$. Since $E^{21}\in V^{\Phi}_{\coh(X)}(2)$ we have $E^{13} = 0$ and so
$E \in V^{\Phi}_{\coh(X)}(2)$.
Therefore $\ch(E^2) = (0,0,0,-r)$.
But it is not possible to have $-r > 0$ and
this is the required contradiction to complete the proof.
\end{proof}
\begin{rem}
Theorem \ref{prop5.1} can be interpreted as saying that if a vector bundle
$E$ over $X$ satisfies $c_1(E)=0=c_2(E)$ then it cannot
carry a non-flat Hermitian-Einstein connection. This is analogous to the case where there are no
charge $1$ $\operatorname{SU}(r)$ instantons on an abelian surface. This is proved in
a slick way using the Fourier-Mukai transform and it would be good to
avoid the direct proof given for Theorem \ref{prop5.1} as it would
follow more directly from Theorem \ref{prop6.9}.
\end{rem}
\section{Auto-equivalences of $\A_{\frac{\sqrt{3}}{2}\ell , \frac{1}{2} \ell}$ under the FMTs}
\label{section6}
Let denote the FMTs $\Psi = L \Phi$ and $\HPsi = \Phi L^{-1} [1]$.
Then by Theorem \ref{prop4.17}, we have that the images of an object
from $\B$ under $\Psi$ and $\HPsi$ are  complexes whose
$\B$-cohomologies can only be non-zero in the $0$, $1$ or $2$ positions.
We  have $\Psi \circ \HPsi \cong (-1)^* \id_{D^b(X)}[-2]$ and $\HPsi
\circ \Psi \cong (-1)^* \id_{D^b(X)}[-2]$. This gives us the following
convergence of spectral sequences.
\begin{specseq}
\label{Bss}
\begin{align*}
E_2^{p,q} &= \Psi^p_{\B} \HPsi^q_{\B} (E) \Longrightarrow H^{p+q-2}_{\B} ((-1)^*E), \\
E_2^{p,q} &= \HPsi^p_{\B} \Psi^q_{\B} (E) \Longrightarrow H^{p+q-2}_{\B} ((-1)^*E),
\end{align*}
for $E$. Here $\Psi^i_{\B} (F):= H^{i}_{\B}(\Psi(F))$ and  $\HPsi^i_{\B} (F):= H^{i}_{\B}(\HPsi(F))$.
\end{specseq}
These convergence of the spectral sequences for $E \in \B$ look
similar to the convergence of some spectral sequences in an abelian
surface for coherent sheaves. See \cite{BBR}, \cite{Macio1},
\cite{Yosh1} for further details.

Recall that if $B_1, B_2 \in \B$ then $\Ext^i(B_1, B_2)  =0$ for any $i<0$.
\begin{prop}
\label{prop6.1}
For $E \in D^b(X)$ we have
\begin{align*}
\Im\,Z(\Psi(E))= -  \Im\,Z(E), \text{ and } \Im\,Z(\HPsi(E)) = - \Im\,Z(E).
\end{align*}
\end{prop}
\begin{proof}
Let $\ch(E) = (a_0, a_1, a_2, a_3)$. Then $\Im\,Z(E) = \frac{3 \sqrt{3}}{4}(a_2 - a_1)$. Also we have
$\ch(\Psi(E)) = (*, a_3 -a_2, a_3-2a_2+a_1, *)$ and
$\ch(\HPsi(E)) = (*, a_2 - 2a_1 +a_0 , -a_1 + a_0 , *)$. Then
$\Im\,Z(\Psi(E)) = \Im\,Z(\HPsi(E)) = - \frac{3 \sqrt{3}}{4}(a_2 - a_1)$ as required.
\end{proof}
From Propositions \ref{prop2.4}, \ref{prop3.1} and Theorem \ref{prop5.1}  we make the following
\begin{note}
\label{prop6.2}
Let $E \in \B$. Then we have the following:
\begin{enumerate}[label=(\Roman*)]
\item if $E \in \HN^\nu(-\infty, 0)$,  then $\mu^+(E_{-1}) < 0$;
\item if $E \in \HN^{\nu}(0, +\infty]$, then $\mu^{-}(E_0) >1$; and
\item for tilt stable $E$ with $\nu(E)=0$, we have
\begin{enumerate}[label=(\roman*)]
\item $\mu^+(E_{-1}) \le 0$,  and $\mu^{-}(E_0) \ge 1$,
\item if $\mu(E_{-1}) =0$ then $E_{-1} = \Po_x$ for some $x \in X$, and
\item if $\mu(E_0) =1$ then $E_0^{**} = L \Po_x$ for some $x \in X$.
\end{enumerate}
\end{enumerate}
\end{note}
\begin{prop}
\label{prop6.3}
Let $E \in \T'$. Then we have the following:
\begin{enumerate}[label=(\roman*)]
\item $H^0_{\coh(X)}(\HPsi^2_{\B}(E)) = 0$, and
\item if $\HPsi^2_{\B}(E) \ne 0$ then $\Im\,Z(\HPsi^2_{\B}(E)) > 0$.
\end{enumerate}
\end{prop}
\begin{proof}
\begin{enumerate}[label=(\roman*)]
\item
For any $x \in X$,
\begin{align*}
\Hom ( \HPsi^2_{\B}(E) , \OO_x  ) & \cong  \Hom ( \HPsi^2_{\B}(E) , \HPsi^2_{\B}( L\Po_{-x} ) ) \\
                                        & \cong  \Hom ( \HPsi(E), \HPsi(L\Po_{-x} ) ) \\
                                        & \cong \Hom ( E, L \Po_{-x} ) = 0,
\end{align*}
since $E \in \T'$ and $L\Po_{-x}  \in \F'$. Therefore $H^0_{\coh(X)}(\HPsi^2_{\B}(E)) = 0$ as required.
\item
From (i), we have $\HPsi^2_{\B}(E) \cong A[1]$ for some $0\neq A \in \HN(-\infty, \frac{1}{2}]$.

Consider the convergence of the spectral sequence:
$$
E^{p,q}_2=\HPsi^{p}_{\coh(X)}(H^{q}_{\coh(X)}(E)) \Longrightarrow \HPsi^{p+q}_{\coh(X)}(E)
$$
for $E$.
Let $E_i := H^{i}_{\coh(X)}(E)$. Then by Note \ref{prop6.2}, $E_0
\in \HN(1, +\infty]$ and so by Corollary \ref{prop4.13} and
  Proposition \ref{prop4.4} we have
$$
(L^{-1}E_0)^2, (L^{-1}E_{-1})^3 \in \HN(0,+\infty].
$$
Therefore from the convergence of the above spectral sequence for $E$, we have
$$
A \in  \HN(-\infty, \frac{1}{2}] \cap \HN(0,+\infty] = \HN(0, \frac{1}{2}].
$$
Let $\ch(A)= (a_0, a_1, a_2, a_3)$. Then from the B-G inequalities for all the H-N semistable factors of $A$, we have
$$
\Im\,Z(\HPsi^2_{\B}(E)) =  \Im\,Z(A[1])  = \frac{3 \sqrt{3}}{4}(a_1 - a_2) > 0
$$
as required.
\end{enumerate}
\end{proof}
\begin{prop}
\label{prop6.4}
Let $E \in \F' $. Then we have the following:
\begin{enumerate}[label=(\roman*)]
\item $H^{-1}_{\coh(X)}(\HPsi^0_{\B}(E)) = 0$, and
\item if $\HPsi^0_{\B}(E) \ne 0$ then $\Im\,Z(\HPsi^0_{\B}(E)) < 0$.
\end{enumerate}
\end{prop}
\begin{proof}
\begin{enumerate}[label=(\roman*)]
\item
Let $x \in X$. Then
\begin{align*}
\Hom ( \HPsi^0_{\B}(E) , \OO_x[1] ) & \cong  \Hom (\Psi \HPsi^0(E) , \Psi (\OO_x[1]))\\
                                    & \cong  \Hom ( \Psi^2_{\B} \HPsi^0_{\B}(E)[-2] , L \Po_x[1] ) \\
                                    & \cong  \Hom ( \Psi^2_{\B} \HPsi^0_{\B}(E) , L \Po_x[3]  ) \\
                                    & \cong  \Hom ( L \Po_x ,  \Psi^2_{\B} \HPsi^0_{\B}(E) )^* .
\end{align*}
From the convergence of the Spectral Sequence \ref{Bss} for $E$, we have the $\B$-SES
$$
0 \to \Psi^0_{\B} \HPsi^1_{\B}(E) \to \Psi^2_{\B} \HPsi^0_{\B}(E) \to F \to 0,
$$
where $F$ is a subobject of $(-1)^* E$ and so $F \in \F'$.
Moreover by the H-N filtration, $F$ fits into the following $\B$-SES
$$
0 \to F_0 \to F \to F_1 \to 0,
$$
where $F_0 \in \HN^{\nu}[0]$ and $F_1 \in \HN^{\nu}( - \infty, 0)$. Since  $L \Po_x \in \HN^{\nu}[0]$,
$$
\Hom (L \Po_x , F_1) =0.
$$
Moreover $F_0$ has a filtration of $\nu$-stable objects $F_{0,i}$ with $\nu (F_{0,i}) = 0$. By Proposition \ref{prop2.6}, each $F_{0,i}$ fits into a non-splitting $\B$-SES
$$
0 \to F_{0,i} \to M_i \to T_i \to 0,
$$
for some $T_i \in \coh^0(X)$ such that $M_i[1] \in \A$ is a minimal object.
Moreover $L \Po_x[1] \in \A$ is a  minimal object.  So finitely many $x \in X$ we can have $L \Po_x \cong M_i$ for some $i$. So for a generic $x \in X$,
$\Hom(L \Po_x, M_i) = 0$ and so $\Hom(L \Po_x , F_{0,i}) = 0$ which implies $\Hom(L \Po_x , F_0) = 0$.
Therefore for a generic $x \in X$, $\Hom ( L \Po_x ,  F ) = 0$.

On the other hand
\begin{align*}
\Hom ( L \Po_x , \Psi^0_{\B} \HPsi^1_{\B} (E)) & \cong  \Hom ( \Psi^0_{\B} (\OO_x)  , \Psi^0_{\B} \HPsi^1_{\B}(E)  )\\
                                               & \cong  \Hom ( \Psi (\OO_x) , \Psi \HPsi^1_{\B}(E) ) \\
                                               & \cong  \Hom (\OO_x , \HPsi^1_{\B}(E) ).
\end{align*}
Here $\HPsi^1_{\B}(E)$ fits into the $\B$-SES
$$
0 \to H^{-1}_{\coh(X)}(\HPsi^1_{\B}(E))[1] \to \HPsi^1_{\B}(E) \to H^0_{\coh(X)}(\HPsi^1_{\B}(E)) \to 0,
$$
where $H^{-1}_{\coh(X)}(\HPsi^1_{\B}(E))$ is torsion free and
$H^0_{\coh(X)}(\HPsi^1_{\B}(E))$ can have torsion supported
on a 0-scheme of finite length. Hence for generic $x \in X$,  $ \Hom (\OO_x , \HPsi^1_{\B}(E) ) = 0$.
Therefore for generic $x \in X$,  $\Hom ( L \Po_x , \Psi^0_{\B} \HPsi^1_{\B} (E))=\Hom( L \Po_x ,  F ) =0$ implies
$\Hom ( L \Po_x ,  \Psi^2_{\B} \HPsi^0_{\B}(E)) =0$. Hence
$\Hom ( \HPsi^0_{\B}(E) , \OO_x[1]  ) = 0$ for generic $x \in X$.
But $H^{-1}_{\coh(X)}(\HPsi^0_{\B}(E))$ is torsion free and so $H^{-1}_{\coh(X)}(\HPsi^0_{\B}(E)) = 0$ as required.

\item
From (i), we have $\HPsi^0_{\B}(E) \cong A$ for some coherent sheaf $0 \ne A \in \HN(\frac{1}{2} , + \infty]$. For any $x \in X$ we have
\begin{align*}
\Ext^1(\OO_x, A) & \cong  \Ext^1 (\OO_x, \HPsi^0_{\B}(E)) \cong \Hom(\Psi(\OO_x), \Psi \HPsi^0_{\B}(E) [1]) \\
                 & \cong  \Hom(L \Po_x,  \Psi^2 \HPsi^0_{\B}(E) [-1]) = 0.
\end{align*}
So $A \in \coh^{\ge 2}(X)$, and if $\ch(A) = (a_0, a_1, a_2, a_3)$ then we have $a_1 > 0$.

Apply the FMT $\Psi$ to $\HPsi^0_{\B}(E)$. Since $\HPsi^0_{\B}(E) \in V^{\Psi}_{\B}(2)$,  $\Psi^2_{\B} \HPsi^0_{\B}(E) \in \B$ has $\coh(X)$-cohomologies:
\begin{itemize}
\item $L A^1$  in position $-1$, and
\item $L A^2$  in position $0$.
\end{itemize}
So we have $A \in V^{\Phi}_{\coh(X)}(1,2)$, $LA^1 \in \HN(-\infty,
\frac{1}{2}]$ and by Corollary  \ref{prop4.13} $A^2 \in \HN(0, +\infty]$.
Therefore  $\ell^2 \ch_1 (A^1) \le 0$ and $\ell^2 \ch_1 (A^2) \ge 0$.
Hence
$$
a_2 \ell^3 = 2 \ell \ch_2 (A) = - \ell^2 \ch_1 (\Phi(A)) = - \ell^2 \ch_1 (A^2) + \ell^2 \ch_1 (A^1) \le 0.
$$
So
$$
\Im\,Z(\HPsi^0_{\B}(E)) =  \Im\,Z(A)  = \frac{3 \sqrt{3}}{4}(a_2 - a_1) < 0
$$
as required.
\end{enumerate}
\end{proof}
\begin{prop}
\label{prop6.5}
\begin{enumerate}[label=(\Roman*)]
\item
Let $E \in \T' $. Then we have the following:\\
(i) $H^{0}_{\coh(X)}(\Psi^2_{\B}(E)) = 0$, and
(ii) if $\Psi^2_{\B}(E) \ne 0$ then $\Im\,Z(\Psi^2_{\B}(E)) > 0$.
\item
Let $E \in \F' $. Then we have the following:\\
(i) $H^{-1}_{\coh(X)}(\Psi^0_{\B}(E)) = 0$, and
(ii) if $\Psi^0_{\B}(E) \ne 0$ then $\Im\,Z(\Psi^0_{\B}(E)) < 0$.
\end{enumerate}
\end{prop}
\begin{proof}
\begin{enumerate}[label=(\Roman*)]
\item Let $E \in \T'$.
\begin{enumerate}[label=(\roman*)]
\item Similar to the proof of (i) in Proposition \ref{prop6.3}.
\item
From (i), we have $\Psi^2_{\B}(E) \cong A[1]$ for some coherent sheaf $0 \ne A \in \HN(-\infty, \frac{1}{2}]$.
Let $\ch(A) = (a_0, a_1, a_2, a_3)$. Then $\ch(L^{-1}A) = (a_0, a_1-a_0, a_2-2a_1+a_0, *)$ and so $a_1 -a_0 <0$.

Apply the FMT $\HPsi$ to $\Psi^2_{\B}(E)$. Since $\Psi^2_{\B}(E) \in V^{\HPsi}_{\B}(0)$,  $\HPsi^0_{\B} \Psi^2_{\B}(E) \in \B$ has $\coh(X)$-cohomologies:
\begin{itemize}
\item $(L^{-1}A)^1$  in position $-1$, and
\item $(L^{-1}A)^2$  in position $0$.
\end{itemize}
So we have  $(L^{-1}A)^2 \in \HN(\frac{1}{2}, +\infty]$, and by Proposition \ref{prop4.12}, $(L^{-1}A)^1 \in \HN(-\infty, 0]$.
Therefore  $\ell^2 \ch_1 ((L^{-1}A)^1) \le 0$ and $\ell^2 \ch_1 ((L^{-1}A)^2) \ge 0$.
Hence
\begin{align*}
(a_2-2a_1+a_0) \ell^3 & =  2 \ell \ch_2 (L^{-1}A) = - \ell^2 \ch_1 (\Phi(L^{-1}A)) \\
                      & =  - \ell^2 \ch_1 ((L^{-1}A)^2) + \ell^2 \ch_1 ((L^{-1}A)^1) \le 0.
\end{align*}
So we have
\begin{align*}
\Im\,Z(\Psi^2_{\B}(E)) & =  \Im\,Z(A[1]) = \frac{3 \sqrt{3}}{4}(a_1- a_2) \\
                        & = - \frac{3 \sqrt{3}}{4} ( (a_1- a_0) + (a_2-2a_1+a_0)) > 0
\end{align*}
as required.
\end{enumerate}
\item Let $E \in \F'$.
\begin{enumerate}[label=(\roman*)]
\item Similar to the proof of (i) in Proposition \ref{prop6.4}.
\item
From (i), we have $\Psi^0_{\B}(E) \cong A$ for some $0\neq A \in \HN( \frac{1}{2}, +\infty]$.

Consider the convergence of the spectral sequence:
$$
E^{p,q}_2=\Psi^{p}_{\coh(X)}(H^{q}_{\coh(X)}(E)) \Longrightarrow \Psi^{p+q}_{\coh(X)}(E)
$$
for $E$.
Let $E_i := H^{i}_{\coh(X)}(E)$. Then by Note \ref{prop6.2},
$E_{-1} \in \HN(-\infty, 0]$ and so by Proposition \ref{prop4.12}
  and Proposition \ref{prop4.4} we have
$$
LE_{-1}^1 \in \HN(-\infty,1], \text{ and } LE_0^0 \in \HN(-\infty,1].
$$
Therefore from the convergence of the above spectral sequence for $E$, we have
$$
A \in  \HN( \frac{1}{2}, +\infty] \cap \HN(-\infty,1] = \HN(\frac{1}{2}, 1].
$$
Also $A$ is reflexive, as $LE_0^0$ and  $LE_{-1}^1$ are reflexive sheaves by Propositions  \ref{prop4.2} and \ref{prop4.6}.
Let $\ch(A)= (a_0, a_1, a_2, a_3)$. Then from the B-G inequalities for all the H-N semistable factors of $A$, we have
$$
\Im\,Z(\Psi^0_{\B}(E)) =  \Im\,Z(A)  = \frac{3 \sqrt{3}}{4}(a_2 - a_1) \le 0.
$$
Equality holds when $A \in \HN[1]$ with $\ch(A) = (a_0, a_0,
a_0, *)$. Then, by considering a Jordan{-}H\"{o}lder filtration
for $A$ together with Theorem \ref{prop5.1}, $L^{-1}A$ has a
filtration of ideal sheaves $\Po_{x_i}\I_{Z_i}$  of some 0-subschemes.
Here $\Psi^0_{\B}(E) \cong A \in V^{\HPsi}_{\B}(2)$ implies $L^{-1}A
\in V^{\Phi}_{\coh(X)}(2,3)$.  An easy induction on the rank of $A$ also shows that
$L^{-1}A\in V^{\Phi}_{\coh(X)}(1,3)$ and so $L^{-1}A\in
V^{\Phi}_{\coh(X)}(3)$. But then  we have $Z_i =\emptyset$ for all $i$.
Therefore $\HPsi^2_{\B}\Psi^0_{\B}(E) \in \coh^0(X)$. Now consider the
convergence of the Spectral Sequence \ref{Bss} for $E$. Then we have
$\B$-SES
$$
0 \to \HPsi^0_{\B}\Psi^1_{\B}(E)  \to \HPsi^2_{\B}\Psi^0_{\B}(E)  \to F \to 0,
$$
where $F$ is a subobject of $(-1)^*E$ and so $F \in \F'$. Then
$\HPsi^2_{\B}\Psi^0_{\B}(E) \in \coh^0(X) \subset \T'$ which implies $F =0$ and
$\HPsi^0_{\B}\Psi^1_{\B}(E)  \cong \HPsi^2_{\B}\Psi^0_{\B}(E)$. But
then we have $\Psi^0_{\B}(E) \cong (-1)^*\Psi^0_{\B} \HPsi^0_{\B}\Psi^1_{\B}(E)  =0$.
This is not possible as  $\Psi^0_{\B}(E) \ne 0$.
Therefore we have the strict inequality $\Im\,Z(\Psi^0_{\B}(E)) <0$ as required to complete the proof.
\end{enumerate}
\end{enumerate}
\end{proof}

\begin{lem}
\label{prop6.6}
\begin{enumerate}[label=(\Roman*)]
\item Let $E \in \T' $. Then (i) $\HPsi^2_{\B}(E) =0$, and (ii) $\Psi^2_{\B}(E) =0$.
\item Let $E \in \F' $. Then (i) $\HPsi^0_{\B}(E) =0$, and (ii) $\Psi^0_{\B}(E) =0$.
\end{enumerate}
\end{lem}
\begin{proof}
\begin{enumerate}[label=(\Roman*)]
\item Let $E \in \T'$.
\begin{enumerate}[label=(\roman*)]
\item From the convergence of the Spectral Sequence \ref{Bss} for $E$, we have the $\B$-SES
$$
0 \to Q \to \Psi^0_{\B} \HPsi^2_{\B}(E) \to \Psi^2_{\B} \HPsi^1_{\B}(E) \to 0.
$$
Here  $Q$ is a quotient of $(-1)^* E \in \T'$ and so $Q \in \T'$.
Then $\Psi^0_{\B} \HPsi^2_{\B}(E)$ fits into the $\B$-SES
$$
0 \to T \to \Psi^0_{\B} \HPsi^2_{\B}(E) \to F \to 0
$$
for some $T \in \T'$ and $F \in \F'$.
Now apply the FMT $\HPsi$ and consider the LES of $\B$-cohomologies.
Then we have $\HPsi^0_{\B}(T) =0$, $\HPsi^1_{\B}(T) \cong \HPsi^0_{\B}(F)$. By  Proposition \ref{prop6.4}
$\Im\,Z(\HPsi^0_{\B}(F)) \le 0$ and by Proposition \ref{prop6.3}  $\Im\,Z(\HPsi^2_{\B}(T)) \ge 0$.
So $\Im\,Z(\HPsi(T)) \ge 0$ and by  Proposition \ref{prop6.1} $\Im\,Z(T)
\le 0$. Since $T \in \T'$, we have $\Im\,Z(T) = 0$ and $\om^2 \ch_1^B(T)=0$. Then by Lemma \ref{prop1.1}, $T \cong T_0$ for some $T_0 \in \coh^0(X)$.
But $\coh^0(X) \subset V^{\HPsi}_{\B}(0)$. Hence $T =0$ and so $Q =0$.
Then $\Psi^0_{\B} \HPsi^2_{\B}(E) \cong \Psi^2_{\B} \HPsi^1_{\B}(E) $ and so we have $\HPsi^2_{\B}(E) \cong (-1)^* \HPsi^2_{\B} \Psi^2_{\B} \HPsi^1_{\B}(E) = 0$ as required.
\item Similar to the proof of (i).
\end{enumerate}
\item Similar to the proofs in (I).
\end{enumerate}
\end{proof}
\begin{cor}
\label{prop6.7}
Let $E \in \B$. Then
\begin{enumerate}[label=(\roman*)]
\item $\Psi^2_{\B}(E), \HPsi^2_{\B}(E) \in \T'$, and
\item $\Psi^0_{\B}(E), \HPsi^0_{\B}(E) \in \F'$.
\end{enumerate}
\end{cor}
\begin{proof}
\begin{enumerate}[label=(\roman*)]
\item
By the definition of $\T'$ and $\F'$,  $\Psi^2_{\B}(E)$ fits into  $\B$-SES
$$
0 \to T \to \Psi^2_{\B}(E) \to F \to 0,
$$
for some $T \in \T'$ and $F \in \F'$. Now apply the FMT $\HPsi$ and consider the LES of $\B$-cohomologies.
Then by Lemma \ref{prop6.6}, $F=0$ as required.

Similarly one can prove $\HPsi^2_{\B}(E) \in \T'$.

\item Similar to the proofs in (i).
\end{enumerate}
\end{proof}
\begin{prop}
\label{prop6.8}
\begin{enumerate}[label=(\Roman*)]
\item Let $E \in \F'$. Then (i) $\HPsi^1_{\B}(E) \in \F'$, and (ii) $\Psi^1_{\B}(E) \in \F'$.
\item Let $E \in \T'$. Then (i) $\HPsi^1_{\B}(E) \in \T'$, and (ii) $\Psi^1_{\B}(E) \in \T'$.
\end{enumerate}
\end{prop}
\begin{proof}
\begin{enumerate}[label=(\Roman*)]
\item
\begin{enumerate}[label=(\roman*)]
\item By the torsion theory $\HPsi^1_{\B}(E)$ fits into   $\B$-SES
$$
0 \to T \to \HPsi^1_{\B}(E) \to F \to 0
$$
for some $T \in \T'$ and $F \in \F'$.
Now we need to show $T=0$.
Apply the FMT $\Psi$ and consider the LES of $\B$-cohomologies.
We get $\Psi^1_{\B}(T) \hookrightarrow \Psi^1_{\B} \HPsi^1_{\B}(E)$ and $T \in V^{\Psi}_{\B}(1)$. Also by the convergence of the Spectral Sequence \ref{Bss} for $E$, $\Psi^1_{\B} \HPsi^1_{\B}(E)$ is a subobject of $(-1)^* E$. Hence $\Psi^1_{\B}(T) \in \F'$ implies $\Im\,Z(\Psi^1_{\B}(T)) \le 0$.
On the other hand by Proposition \ref{prop6.1},
$ \Im\,Z(\Psi^1_{\B}(T)) = \Im\,Z(T) \ge 0$ as $T \in \T'$. Hence $\Im\,Z(T) = 0$ and  $T \in \T'$ implies $\om^2 \ch_1^B (T)  = 0$.
So by Lemma \ref{prop1.1}, $T \cong T_0$ for some $T_0 \in \coh^0(X)$. Since any object from $\coh^0(X)$ belongs to $V^{\Psi}_{\B}(0)$, $\Psi^1_{\B} (T) = 0$. So $T=0$ as required.
\item Similar to the proof of (i).
\end{enumerate}
\item Similar to the proofs in (I).
\end{enumerate}
\end{proof}

By Lemma \ref{prop6.6}, Corollary \ref{prop6.7} and Proposition \ref{prop6.8} we have
$$
\Psi[1]\left( \F'[1] \right) \subset \A, \ \text{ and }  \ \Psi[1]\left( \T' \right) \subset \A.
$$
Since $\A = \langle \F'[1], \T' \rangle$,  $\Psi[1]\left( \A \right) \subset \A$.

Similarly we have $\HPsi[1]\left( \A \right) \subset \A $.
The isomorphisms $\HPsi[1] \circ \Psi[1] \cong (-1)^* \id_{D^b(X)}$ and $\Psi[1] \circ \HPsi[1] \cong (-1)^* \id_{D^b(X)}$ give us the following

\begin{thm}
\label{prop6.9}
The FMTs $\Psi[1]$ and $\HPsi[1]$ give the auto-equivalences
$$
\Psi[1]\left( \A \right) \cong \A, \text{ and } \HPsi[1]\left( \A \right) \cong \A
$$
of the abelian category $\A$.
\end{thm}
\section*{Acknowledgements}
The authors would like to thank Arend Bayer and Tom Bridgeland for very useful discussions and comments.
We are also grateful to Jason Lo and  Yukinobu Toda for pointing out several errors.
Special thanks go to the referee for a thorough reading of the paper and insightful comments
that led to a substantial improvement of this paper, especially in sections 4 and 5.
The second author is funded by Principal's Career Development Scholarship programme and Scottish Overseas
Research Student Awards Scheme of the University of Edinburgh, and this work forms a part
of his PhD thesis.




\end{document}